\documentclass[12pt]{article}
 
\usepackage{amssymb} 
\usepackage{a4} 
\newcommand{\sect}[1]{\setcounter{equation}{0}\section{#1}}

\newcommand{\app}[1]{\setcounter{section}{0} 
\setcounter{equation}{0} \renewcommand{\thesection} 
{\Alph{section}}\section{#1}} 
 
\newcommand{\A}{{\cal A}} 
 
\newcommand{\la}{\langle} 
\newcommand{\ra}{\rangle} 
\def\1{{\bf 1}} 
\def\b#1{{\mathbb #1}} 
\def\c#1{{\cal #1}} 
\def\u#1{{\underline #1}}
\def\R{{\cal R}\,} 
\def\F{{\cal F}\,} 
 
\newcommand{\tl}{\,\triangleleft} 

\def\cross{{\triangleright\!\!\!<}} 
\def\cocross{{>\!\!\!\triangleleft\,}} 
 
\def\id{\mbox{id\,}}

\newcommand{\bphi}{\mbox{\boldmath $\phi$}}
\newcommand{\bpsi}{\mbox{\boldmath $\psi$}}
\newcommand{\bpi}{\mbox{\boldmath $\pi\!$}}

\newcommand{\cbphi}{\mbox{\boldmath $\check\phi$}} 
\newcommand{\cbpsi}{\mbox{\boldmath $\check\psi$}}

\newcommand{\hbphi}{\mbox{\boldmath $\hat\phi$}}  
\newcommand{\hbpsi}{\mbox{\boldmath $\hat\psi$}} 
\newcommand{\balpha}{\mbox{\boldmath $\alpha$}}   
\newcommand{\bbeta}{\mbox{\boldmath $\beta$}} 
\newcommand{\calpha}{\mbox{\boldmath $\check\alpha$}}  
\newcommand{\cbeta}{\mbox{\boldmath $\check\beta$}}
 
\def\nn{\nonumber \\} 
\newcommand{\be}{\begin{equation}} 
\newcommand{\ee}{\end{equation}} 
\newcommand{\bea}{\begin{eqnarray}} 
\newcommand{\eea}{\end{eqnarray}} 
\newcommand{\ba}{\begin{array}} 
\newcommand{\ea}{\end{array}} 
\newtheorem{prop}{Proposition} 
\newtheorem{lemma}{Lemma} 
\newtheorem{theorem}{Theorem}

\def\sq{\mbox{\rlap{$\sqcap$}$\sqcup$}} 
\newenvironment{proof}[1]{\vspace{5pt}\noindent{\bf Proof #1}\hspace{6pt}}%
{\hfill\sq} 
\newcommand{\bp}{\begin{proof}} 
\newcommand{\ep}{\end{proof}\par\vspace{10pt}\noindent} 

\begin{document}

\title{On the hermiticity of $q$-differential operators and forms
on the quantum Euclidean spaces $\b{R}_q^N$}

\author{        Gaetano Fiore   \\\\
         \and
        Dip. di Matematica e Applicazioni, Fac.  di Ingegneria\\ 
        Universit\`a di Napoli, V. Claudio 21, 80125 Napoli
\\
        and
\\
        I.N.F.N., Sezione di Napoli,\\
        Complesso MSA, V. Cintia, 80126 Napoli
        }

\date{}

\maketitle

\abstract{We show that
the complicated $\star$-structure characterizing for positive $q$
the $U_qso(N)$-covariant differential calculus on the
non-commutative manifold $\b{R}_q^N$
 boils down to similarity transformations
involving the ribbon element
 of a central extension of
$U_qso(N)$ and its formal square root $\tilde v$. 
Subspaces of the spaces of functions and of $p$-forms on $\b{R}_q^N$ are made 
into Hilbert spaces by introducing non-conventional
``weights'' in the integrals defining the corresponding scalar products,
namely suitable positive-definite $q$-pseudodifferential operators
$\tilde v'{}^{\pm 1}$ realizing the action of $\tilde v^{\pm 1}$; 
this serves to
 make the partial $q$-derivatives 
antihermitean and the exterior coderivative equal to the
hermitean conjugate of the exterior derivative, as usual.
There is a residual freedom
in the choice of the weight $m(r)$
 along the `radial coordinate' $r$. Unless we choose a constant
 $m$, then the square-integrables functions/forms mustfulfill an additional
 condition, namely their analytic continuations to thecomplex $r$ plane
 can have poles only on the sites of some special lattice.Among the functions
 naturally selected by this condition there are $q$-special functions with
 `quantized' free parameters.}
 
\vskip2cm

MSC-class: 81R50; 81R60; 16W10; 16W30; 20G42.


\newpage

\sect{Introduction}

Over the past two decades the noncommutative geometry program
\cite{Con} and the related program  of generalizing the concept of
symmetries through quantum  groups \cite{Dri86,Wor87-89,FadResTak89} and
quantum group covariant noncommutative spaces (shortly: quantum spaces)
\cite{Man,FadResTak89} has found a widespread interest in the mathematical 
and theoretical physics  community and 
accomplished substantial progress. %
Initially,  mathematical investigations have been concentrated
essentially in compact noncommutative manifolds, the
non-compact being
usually much more  complicated to deal with, especially when
trying to proceed from an algebraic to a functional-analytical treatment.
 In particular, so are $\star$-structures and 
$\star$-representations  of the involved algebras. 
Recently, an increasing number of works is being devoted to
extend results to non-compact noncommutative manifolds. We might
divide these works into two subgroups.  The first (see e.g.
\cite{ConDub02,GayIocVar05,GayGraIocSchVar04,GraLizMarVit02,Sui04})  essentially
deal with  non-compact noncommutative  manifolds which can be obtained by
isospectral deformations \cite{ConLan01}
of commutative Connes' spectral triples and carry
the action of an abelian group $\b{T}^k\times\b{R}^h$.
The second, and even more difficult  (see  e.g. \cite{MasNakWor03}, and
references therein) deal with  non-compact noncommutative  manifolds which
underlie some quantum group or more generally carry the action of some 
quantum group; it is still under debate what the most convenient
axiomatization  of these models is \cite{MasNakWor03}. 

The noncommutative manifold we are going to consider in the present
work belongs to the second category and is relatively old and famous, but
presents an additional complication even at the formal level (i.e. before
entering a functional-analytic treatment):
the $\star$-structure characterizing for
real $q$  the $U_qso(N)$-covariant differential calculus \cite{CarSchWat91}
on the quantum
Euclidean space $\b{R}_q^N$ \cite{FadResTak89} is characterized 
by an unpleasant 
nonlinear action  on the
differentials, the partial derivatives and the exterior derivative
\cite{Ogi92}. This at the origin of a host
of formal and substantial complications.
As examples we mention the following
difficulties: determinining the actual geometry of $\b{R}_q^N$
\cite{FioMad99,CerFioMad00};
identifying  the `right' momentum sector
within the algebra of observables of quantum mechanics
on a $\b{R}_q^N$-configuration space and  solving the
 corresponding
eigenvalue problems for Hermitean operators in the form of differential
operators \cite{Wei94,Fio95JMP,Wes99};
more generally formulating and
solving differential equations on $\b{R}_q^N$;
finally, writing down tractable kinetic terms
for Lagrangians of potential field theory models on $\b{R}_q^N$.
A similar situation occurs for other non-compact quantum spaces, notably
for the $q$-Minkowski space \cite{OgiSchWesZum92}.

It turns out that we are facing a
problem similar to the one we encounter in functional analysis
on the real line when taking
 the Hermitean conjugate of a differential operator like
\be
D=\sigma(x)\, \frac d{dx}\,\frac 1{\sigma(x)},
   \label{undef}
\ee
where $\sigma(x)$ is a smooth
 complex function vanishing for no $x$. 
As an element of the Heisenberg algebra $D$ 
is not imaginary
 (excluding the trivial case 
$\sigma\equiv 1$)
 w.r.t. the $\star$-structure 
$$
x^{\star}=x, \qquad \left(\frac d{dx}\right)^{\star}=-\frac d{dx},
$$
but fulfills the similarity transformation
$$
D^{\star}=-|\sigma|^{-2} D\,|\sigma|^2,
$$
this corresponding to the fact that it is not
antihermitean as an operator
on $L^2(\b{R})$. $D$ is however (formally) antihermitean
on $L^2(\b{R},|\sigma|^{-2} dx)$.  In other words, if we 
  insert
the weight $|\sigma|^{-2}>0$ in the integral giving the
scalar product, 
$$
(\bphi,\bpsi)=\int \bphi^{\star}(x)\,|\sigma|^{-2}\,\bpsi(x)dx,
$$
[as one does when setting the Sturm-Liouville problem
for  $D^2$], $D$ becomes antihermitean under the corresponding Hermitean
conjugation $\dagger$\footnote{The Hermitean conjugation $\dagger$ is
the representation of the following modified $\star$-structure $\star'$ of the
Heisenberg algebra 
$a^{\star'}=[|\sigma|^{-2}\,a\,
|\sigma|^2]^{\star}=|\sigma|^2\,a^*\,|\sigma|^{-2}$.
 }:
$$
( A^{\dagger}\bphi,\bpsi):=( \bphi,A\bpsi)
\qquad\quad\Rightarrow\qquad\quad D^{\dagger}=-D.
$$

In this work we show that the partial derivatives
$\partial^{\alpha}$ and the exterior derivative $d$ of the 
$U_qso(N)$-covariant differential calculus on  $\b{R}_q^N$
can be expressed by the similarity transformation
\be
\partial^{\alpha}=
\tilde \nu '{}\tilde\partial^{\alpha} 
\tilde \nu '{}^{-1},  \qquad \qquad \qquad d=  \tilde \nu'{}^{-1}{}\tilde 
d\tilde \nu',                       \label{simila}
\ee
in terms of elements $\tilde\partial^{\alpha},\tilde d$ which are 
purely imaginary under the $\star$-structure  studied 
in \cite{Ogi92}. 
The unusual and novel feature here is that 
$\tilde \nu '{}$ is not a function on $\b{R}_q^N$ but a  
positive-definitepseudodifferential operator, more precisely  
the realization of the fourth 
root of the ribbon element of the extension of  
$U_qso(N)$ with a central element generating  
dilatations of $\b{R}_q^N$.  
Therefore the $\partial^{\alpha}$ become  antihermitean  
and the exterior coderivative $\delta$ becomes 
the Hermitean conjugate of $d$ (on the 
space of differential forms) if we introduce the  
``weights'' $\tilde v'{}^{\mp 1}:=\tilde \nu '{}^{\mp 2}$ in the integral
defining the   scalar product of two `wave-functions/forms' on 
$\b{R}_q^N$.  

For practical purposes it is much more 
convenient to use the $\partial^{\alpha}$ rather than the 
$\tilde\partial^{\alpha}$  because the former have
much simpler commutation relations
(in the form of modified Leibniz rules) with the 
coordinates  of $\b{R}_q^N$, 
whereas for the commutation relations involving the 
$\tilde\partial^{\alpha}$  we
have even  not found a closed form. This suggests to cure
the complications mentioned at the beginning as one does in the
undeformed, functional-analytical setting.

Section \ref{preli} contains preliminaries about the quantum group
$U_qso(N)$, the
 differential calculus on $\b{R}_q^N$, frame bases,
Hodge map and the analog of Lebesgue integration 
over $\b{R}_q^N$; the latter is completely determined apart form a residual 
freedom in choosing the integration measure $m(r)dr$ 
along the radial direction $r$. 
In section \ref{simi} we prove at the algebraic level
(i.e. at the level of formal power series) Eq.
(\ref{simila}) and the corresponding formula
for the differentials
$dx^i$ of the coordinates $x^i$ of $\b{R}_q^N$.
In section \ref{hilbert} we deal with implementing the
previous algebraic results in a functional-analytical
setting:  we introduce
spaces of square-integrable functions/forms over $\b{R}_q^N$ 
and show how the algebraic $\star$-structure can be implemented 
in different ``pictures'' (i.e. configuration space realizations) as
Hermitean 
 conjugation of operators acting on them. As applications, we first consider quantum mechanics on
$\b{R}_q^N$ and recall how one can diagonalize  a set of commuting
observables including various momentum components, then we write down
`tractable' kinetic terms for (bosonic) field theories on $\b{R}_q^N$.  
These steps
require promoting the formally (i.e. algebraically) defined 
$\tilde\nu '{}^{\pm 1}$ into corresponding well-defined 
pseudodifferential 
operators, and this is done in section \ref{PDHPO} passing to the
Fourier transform of the variable $y=\ln r$. No further
constraint is needed  if $m(r)\equiv 1$, whereas an additional
one must be imposed on the spaces of square-integrable functions/forms if 
$m(r)$ is not constant (non-homogeneous space along the radial direction),
e.g. if $m(r)dr$ is the measure of the socalled Jackson integral:
they have to be restricted
to interesting subspaces $L_2^m$ consisting of functions whose analytic
continuation in the complex $r$-plane have poles locations $r_{\alpha}$
on a certain number $\gamma$ of ``rays'' originating from $r=0$, forming with
each  other angles equal to $2\pi/\gamma$, 
and such that $|r_{\alpha}|=q^j$ (or 
$|r_{\alpha}|=q^{j+\frac 12}$),  with $j\in\b{Z}$. Surprisingly,
this is a condition which automatically selects $q$-special functions
where their free parameters (which will play the role of fundamental
physical quantities, e.g. a universal energy scale) are
``quantized''.

\sect{Preliminaries}

\label{preli}

\subsection{$\b{R}_q^N$ and its covariant 
differential calculi}

\label{preli1} 

As a noncommutative space we consider the $U_qso(N)$-covariant deformation
\cite{FadResTak89} of the Euclidean space
$\b{R}^N$ ($h:=\ln q$ plays the role of deformation parameter). 
We shall call the deformed algebra 
of functions on this space  ``algebra of functions on the quantum
Euclidean space $\b{R}_q^N$'', and denote it by $F$. 
It is essentially the unital associative
algebra over $\b{C}[[h]]$ generated by $N$ elements $x^i$ (the cartesian
``coordinates'') modulo the relations 
(\ref{xxrel}) given below, and
will be extended to include formal power series in the generators;
out of $F$ we shall extract subspaces consisting of elements
that can be considered integrable or square-integrable functions.
The $U_qso(N)$-covariant differential calculus
on $\b{R}_q^N$ \cite{CarSchWat91} is defined introducing 
the invariant exterior derivative $d$, satisfying nilpotency
and the Leibniz rule $d(fg)=dfg+fdg$, and imposing the covariant
commutation relations
 (\ref{xxirel}) between the $x^i$ and the differentials
$\xi^i:=dx^i$. 
 Partial derivatives are introduced through the 
decomposition
$d=:\xi^i\partial_i$. All the other commutation relations
are derived by consistency. The complete list is
\bea
&& {\cal P}_a{}^{ij}_{hk}x^hx^k=0, \label{xxrel}\\
&& x^h\xi^i=q\hat R^{hi}_{jk}\xi^jx^k,\label{xxirel}\\
&& ({\cal P}_s+{\cal P}_t)^{ij}_{hk}\xi^h\xi^k=0,\label{xixirel}\\
&& {\cal P}_a{}^{ij}_{hk}\partial_j\partial_i=0, \label{ddrel}\\
&& \partial_i x^j = \delta^j_i+q\hat R^{jh}_{ik} x^k\partial_h, 
                                                 \label{dxrel}\\
&& \partial^h\xi^i=q^{-1}\hat R^{hi}_{jk}\xi^j \partial^k.\label{dxirel}
\eea
The $N^2 \times N^2$ matrix $\hat{R}$ is the braid matrix of 
$SO_q(N)$ \cite{FadResTak89}.
The matrices
${\cal P}_s$, ${\cal P}_a$, ${\cal P}_t$ are $SO_q(N)$-covariant
deformations of the symmetric trace-free,
antisymmetric and trace projectors respectively, which appear in the
projector decomposition of $\hat{R}$
\be
\hat R = q{\cal P}_s - q^{-1}{\cal P}_a + q^{1-N}{\cal P}_t.       
                                                \label{projectorR}
\ee
%
The ${\cal P}_t$ projects on a
one-dimensional sub-space and 
 can be written in the form 
\be
{\cal P}_t{}_{kl}^{ij} = (g^{sm}g_{sm})^{-1} g^{ij}g_{kl}
= \frac{q^2-1}{(q^N-1)(1+q^{2-N})} g^{ij}g_{kl}  \label{Pt}
\ee
where the $N \times N$ matrix $g_{ij}$ 
is a $SO_q(N)$-isotropic tensor, 
deformation of the ordinary Euclidean metric.
The metric and the braid matrix
satisfy the 
 relations~\cite{FadResTak89}
\be
g_{il}\,\hat R^{\pm 1}{}^{lh}_{jk} = 
\hat R^{\mp 1}{}^{hl}_{ij}\,g_{lk}, \qquad
g^{il}\,\hat R^{\pm 1}{}_{lh}^{jk} = 
\hat R^{\mp 1}{}_{hl}^{ij}\,g^{lk}.                          \label{gRrel}
\ee
Indices will be lowered and raised using $g_{ij}$ and its inverse $g^{ij}$,
e.g.
$$
\partial^i:=g^{ij}\partial_j \qquad \qquad
x_i:=g_{ij}x^j.
$$

\medskip

We shall call ${\cal DC}^*$ (differential calculus algebra
on $\b{R}_q^N$) the unital associative 
algebra  over $\b{C}[[h]]$ generated by 
$x^i,\xi^i,\partial_i$ modulo these relations.
We shall denote by 
$\bigwedge^*$ (exterior algebra, or algebra of exterior forms)
the graded unital subalgebra generated by the $\xi^i$ alone, with
grading $\natural \equiv$the degree in $\xi^i$, and by $\bigwedge^p$ 
(vector
space of exterior $p$-forms) the component with grading $\natural =p$, 
$p=0,1,2...$. 
Each $\bigwedge^p$ carries an irreducible representation
of $U_qso(N)$, and its dimension is the binomial
coefficient $N\choose{p}$ \cite{Fio94}, exactly as in the $q=1$
(i.e. undeformed) case; in particular there are no forms with
$p>N$, and $\mbox{dim}(\bigwedge^N)={N\choose{N}}=1$, therefore
$\bigwedge^N$ carries the singlet representation of $U_qso(N)$. 

We shall endow ${\cal DC}^*$ with the same grading $\natural $, and call
${\cal DC}^p$ its component with grading $\natural =p$. The elements of
${\cal DC}^p$ can be considered differential-operator-valued $p$-forms.

We shall denote by $\Omega^*$ (algebra of differential forms)
the graded unital subalgebra generated by the $\xi^i,x^i$, with
grading $\natural $, and by $\Omega^p$ 
(space of differential $p$-forms) its component with grading $p$; 
by definition $\Omega^0=F$ itself.
Clearly both $\Omega^*$ and $\Omega^p$ are
$F$-bimodules.

We shall denote by ${\cal H}$ (Heisenberg algebra on $\b{R}_q^N$)
the unital subalgebra generated by the $x^i, \partial_i$. Note that
by definition ${\cal DC}^0={\cal H}$, and that both
${\cal DC}^*$ and ${\cal DC}^p$ are ${\cal H}$-bimodules.

Using (\ref{ddrel}), (\ref{gRrel}) one can easily 
verify that the $\partial^i$
satisfy the same commutation relations as the $x^i$, and therefore
together with the unit $\1$ generate 
a subalgebra of ${\cal H}$ isomorphic to $F$, 
which we shall call $F'$.
Denote by 
$\{{\cal D}_{\pi}\}_{\pi\in \Pi}$ a basis of the vector space 
underlying $F'$
consisting of homogeneous polynomials in the $\partial$'s
and with first element ${\cal D}_0=\1$.
Any ``pseudodifferential-operator-valued form'', i.e. any element 
${\cal O}\in {\cal DC}^*$, 
(in particular  ${\cal O}\in {\cal H}$) 
can be uniquely expressed in the ``normal-ordered'' form
\be
{\cal O}=\sum_{\pi\in \Pi}{\cal O}_{\pi}{\cal D}_{\pi},
\qquad \qquad {\cal O}_{\pi}\in \Omega^*           \label{Odeco}
\ee
by repeated application of 
relation (\ref{dxrel}), (\ref{dxirel}) to move step by step all $\partial$'s
to the right of all $x,\xi$'s. For any $\omega\in \Omega^*$ 
we shall
denote by ${\cal O}\omega|$ the $\pi=0$ component $({\cal O}\omega)_0$ 
of the  normal-ordered form of ${\cal O}\omega$:
$$
{\cal O}\omega=\sum_{\nu\in \Pi}({\cal O}\omega)_{\nu}{\cal D}_{\nu}
={\cal O}\omega|+\sum_{\nu\neq 0}({\cal O}\omega)_{\nu}{\cal D}_{\nu}.
$$
In particular, for 
${\cal O}=\partial_i$ and $\omega\equiv f\in F$ the previous
formula becomes
 the deformed Leibniz rule
\be
\partial_if= \partial_if|+ f_i^j\partial_j,\qquad\qquad f_i^j\in F.
   \label{Leibniz}
\ee
{}From (\ref{dxrel}) we find e.g. that if $f=x^h$ then
$\partial_if|=\delta^h_i$ and $f_i^j=q\hat R^{hj}_{ik}x^k$.
We have introduced this vertical bar $|$ in the
notation
 to make always clear ``where the action of the derivatives
is meant to stop'', while sometimes this remains ambiguous
by the mere use of brackets. From associativity the
obvious property
$$
{\cal O}({\cal O}'\omega|)|={\cal O}{\cal O}'\omega|
$$
follows. $F,F'$ are 
dual vector spaces w.r.t. the pairing  \cite{Maj95}
\be
\la \partial_{i_1}...\partial_{i_l},
x^{j_1}...x^{j_m}\ra =\delta_{lm}  \partial_{i_1}...\partial_{i_l}
x^{j_1}...x^{j_l}|  \in\b{C}       \label{pair}
\ee
with $m=0,1,...$.

\medskip

The elements 
$$
r^2\equiv x \cdot x:=x^k x_k, \qquad \qquad 
\partial\cdot\partial:=
g^{kl} \partial_l \partial_k=\partial^k \partial_k
$$ 
are $U_qso(N)$-invariant and respectively generate the centers
of $F,F'$. $\partial\cdot\partial$ is a deformation of
the Laplacian on $\b{R}^N$. 
We shall slightly extend $F$ by
introducing the square root $r$ of $r^2$ 
and its inverse $r^{-1}$ as new (central) generators;
$r$ can be considered as the deformed ``Euclidean
distance  of the generic point of coordinates $(x^i)$ of $\b{R}_q^N$ 
from the origin''.
Then the elements $t^i:=x^i r^{-1}$ fulfill (\ref{xxrel}) as well as
the relation $t\cdot t=1$; they generate the deformed algebra 
$F(S_q^{N\!-\!1})$ of
``functions on the unit quantum Euclidean sphere''. The latter can be
completely decomposed into eigenspaces $V_l$ of the deformed quadratic Casimir of
$U_qso(N)$,  or equivalently of the Casimir $w$ defined in (\ref{defw})
with eigenvalues $w_l:=q^{-l(l+N-2)}$, implying a corresponding 
decomposition for 
$F$:
\be
F(S_q^{N\!-\!1})=\bigoplus\limits_{l=0}^{\infty}V_l \qquad\qquad
F=\bigoplus\limits_{l=0}^{\infty}
\left(V_l\otimes \b{C}[[r,r^{-1}]]\right) \label{ldeco} 
\ee

An orthonormal basis $\{S_l^I\}$ (consisting of
`spherical harmonics') of $V_l$ 
can be extracted out of the set of homogeneous,
completely symmetric and trace-free polynomials
of degree $l$
\be
S_l^I\equiv S_l^{i_1i_2...i_l}:={\cal P}^{s,l}{}^{i_1i_2...i_l}_{j_1j_2...j_l}
t^{j_1}t^{j_2}...t^{j_l}                \label{defS}
\ee
suitably normalized ($I$ denotes the multi-index $i_1i_2...i_l$, 
${\cal P}^{s,l}$ denotes the $U_qso(N)$-covariant,
completely symmetric and trace-free projector with $l$ indices
\cite{Fio93,Fio04}).
 Therefore
for the generic $f\in F$
\be
f=\sum\limits_{l=0}^{\infty} f_l=\sum\limits_{l=0}^{\infty}
\sum_IS_l^I f_{l,I}(r).                       \label{ldeco2} 
\ee

\bigskip

The ${\star}$-structure compatible with the compact
${\star}$-structure of $U_qso(N)$ requires 
$q\in \b{R} \setminus\{0\}$.  On the generators $x^i$ 
${\star}$ is 
given by \cite{FadResTak89}\footnote{If we enumerate the $x^i$ of \cite{FadResTak89}  as in \cite{Ogi92}  
by $i=-n,\ldots,-1,0$, $1,\ldots n$ for $N$ odd,  
and $i=-n,\ldots,-1, 1,\ldots n$ for $N$ even, where   
 $n:=\left[\frac N 2\right]$ is the rank of   
$so(N)$, then  the metric matrix 
reads $g_{ij}=g^{ij}=q^{-\rho_i}\delta_{i,-j}$,  where 
$(\rho_i):=\Big(\frac N 2\!-\!1,\frac N 2\!-\!2,\ldots,\frac{1}{2},0, 
\!-\!\frac{1}{2}, \ldots,1\!-\!\frac N 2\Big)$  
for $N$ odd,   
$(\rho_i):=\Big(\frac N 2\!-\!1,\frac N 2\!-\!2,\ldots,0,0,\ldots,1\!-\!\frac N 2\Big)$ 
for $N$ even.  We can obtain a set of $N$ real coordinates   $x^{\alpha}$ by
 a linear transformation $x^{\alpha}:=V^{\alpha}_ix^i$ 
($\alpha=0,1,...,2n$ for odd $N$, $\alpha=1,...,2n$   
for even $N$) defined by ($h\ge 1$) 
$$ 
V^{2h-1}_i\!:= \!\frac 1{\sqrt{2}}(\delta^h_i\!+\!g_{ih}),\qquad\qquad  
V^{2h}_i\!:= \!\frac {-i}{\sqrt{2}}(\delta^h_i\!-\!g_{ih}), \qquad\qquad  
V^0_i\!:=\!\delta^0_i \qquad\quad\mbox{(only for odd $N$)}. 
$$ 
}  
\be
x^i{}^{\star}=x^jg_{ji}                      \label{starx}         
\ee
whereas the conjugates of the derivatives $\partial^i$
resp. the differentials) are not 
combinations of the derivatives (resp. the differentials) themselves. 
One can complete a $U_qso(N)$-covariant  
${\star}$-structure by the relations \cite{OgiZum92} 
\be
\xi^i{}^{\star}=\hat\xi^jg_{ji}    \qquad\qquad
\partial^i{}^{\star}=-q^{-N}\hat\partial^j g_{ji}, \label{starxid}                    
\ee
where 
\bea
\hat\partial^i &:=& \Lambda^2\left[\partial^i+
\frac{q\,k}{1+q^{2-N}}x^i\partial\cdot\partial\right],\qquad\qquad k:=q\!-\!q^{-1}\\
\hat\xi^i &:=& \sigma q^N\Lambda^{-2}\left[\xi^i+q^{-1}kx^id-k
\left(q^{1-N}\xi\cdot x+\frac{k\,q^{-2}}{1+q^{N-2}}r^2\, d\right)
\hat\partial^i\right]                          \label{barredxi}\\
&=& \sigma q^{N-2}\Lambda^{-2}\left[\xi^i+qk\xi^j\partial_jx^i-k
\left(q^{1-N}\xi\cdot x+\frac{k}{1+q^{N-2}}\xi^j\partial_jr^2\right)
\hat\partial^i\right];\nonumber
\eea
the second expression in (\ref{barredxi}) is derived from the
first \cite{OgiZum92} using the Leibniz rule and the decomposition
$d=\xi^i\partial_i$. Here
$\sigma$ is a pure phase factor which we shall set $=1$, 
whereas the element $\Lambda^{-2}$ is defined by 
\be
\Lambda^{-2}:=1 +qkx^i\partial_i+
\frac{q^N\,k^2}{(1+q^{N-2})^2}r^2\partial\cdot\partial\equiv 1+O(h)
\ee
(in 
\cite{OgiZum92} it was denoted by $\Lambda$).
Its square root and inverse square root
$\Lambda^{-1},\Lambda$ can be either introduced as additional
generators or as formal power series in the deformation parameter
$h=\ln q$. They fulfill the relations
\be
\Lambda x^i=q^{-1}x^i\Lambda,\qquad
\Lambda\partial^i=q\partial^i\Lambda, \qquad
\Lambda \xi^i=\xi^i\Lambda,  \qquad \Lambda 1|=1\quad\label{Lambdaprop}
\ee
and the corresponding ones for $\Lambda^{-1}$.
The elements $\hat\xi^i,\hat\partial_i$ 
satisfy relation (\ref{xixirel}-\ref{ddrel}) and the analogue of 
(\ref{dxrel}-\ref{dxirel}) with 
$q,\hat R$ replaced by $q^{-1},\hat R^{-1}$. As a consequence
$\hat d:=\hat\xi^i\hat\partial_i= -d^{\star}$ is also
$U_qso(N)$-invariant, nilpotent, and satisfies the
Leibniz rule on $F$. In fact $\hat d,\hat\xi^i,\hat\partial_i$
can be introduced also as independent objects defining an alternative
$U_qso(N)$-covariant differential calculus. We shall denote by 
$\hat F'$ the subalgebra generated by the $\hat\partial_i$; it is
isomorphic to $F,F'$, too.
One finds \cite{OgiZum92} that under the action of $\star$
\be
r^{\star}=r,\qquad(\partial\cdot\partial)^{\star}=
q^{-2N}\hat\partial\cdot\hat\partial=q^{2-N}\partial\cdot\partial
\Lambda^2, \qquad
\Lambda^{\star}=q^N\Lambda^{-1}.
\quad\ee

\subsection{$\widetilde{U_qso(N)}$ and its action on ${\cal DC}^*$}
\label{vf}

We extend as in Ref. \cite{Maj94} the compact Hopf $\star$-algebra $U_qso(N)$  
(this requires real $q$) by adding a central, primitive and imaginary generator
$\eta$ 
 $$ 
\Delta(\eta)=\1\otimes\eta+\eta\otimes\1, \qquad\qquad 
\epsilon(\eta)=0,\qquad\qquad S\eta=-\eta 
$$ 
(here $\Delta,\epsilon,S$ respectively denote the coproduct, counit, antipode),
and we endow the resulting Hopf $\star$-algebra $H:=\widetilde{U_qso(N)}$  
by the  quasitriangular structure
\be 
\tilde\R:=\R\, q^{\eta\otimes\eta}, 
\ee 
where $\R\equiv\R^{(1)}\otimes\R^{(2)}$ (in a Sweedler notation with upper
indices and suppressed  summation index) denotes the  quasitriangular
structure of $U_qso(N)$. 
This ${\star}$-structure of $H$ thus 
can be summarized by the relations
\be 
\R^{(1)}{}^{\star}\otimes\R^{(2)}{}^{\star}=\R_{21},  
\qquad \qquad \eta^{\star}=-\eta.\label{compact}
\ee 
$\mathcal{DC}^*$ is $H$-module $\star$-algebra (which here 
we choose to be {\it right}),
\be 
(aa')\tl g = (a\!\tl g_{(1)})\, (a'\!\tl g_{(2)}).  \label{modalg2} 
\ee 
Here $g_{(1)} \otimes g_{(2)}=\Delta(g)$ in Sweedler notation.
The transformation laws of the generators
$\sigma^i=x^i,\xi^i,\partial^i$ of
 ${\cal DC}^*$ under the $H$-action read
\bea 
&&\sigma^i\tl \, g=\rho^i_j(g)\sigma^j\qquad\qquad  
g\in U_qso(N), \qquad \label{fundrep}\\ 
&&x^i\tl \,\eta=x^i,\qquad\qquad \xi^i\tl \,\eta=\xi^i, 
\qquad\qquad \partial^i\tl \,\eta=-\partial^i; \qquad 
\eea 
here $\rho$ denotes the $N$-dimensional representation of $U_qso(N)$. 
The braid matrix $\hat R$ is related to $\R$ by  
$\hat R^{ij}_{hk}=\rho^j_h(\R^{(1)})\,\rho^i_k(\R^{(2)})$;
its explicit form can be found in \cite{FadResTak89}. The
elements 
$$
Z^i_j:= T^{(1)}\!\rho^i_j(T^{(2)}), \qquad \quad       
\mbox{where }T\!=\!\R_{21}\R\!\equiv \!T^{(1)}\!\otimes\! T^{(2)}\!,  
\quad\R_{21}\!\equiv\!\R^{(2)}\!\otimes\!\R^{(1)}
$$  
are generators of $U_qso(N)$, and make up the
``$SO_q(N)$ vector field matrix''
$Z$  \cite{Zum91,Zum92,SchWatZum91,SchWatZum92}. 
The $Z^i_j$ are related  to the Faddeev-Reshetikin-Takhtadjan  
generators \cite{FadResTak89} 
\be  
{\cal L}^{+,}{}_l^a:=\R^{(1)}\rho_l^a(\R^{(2)}),\qquad\qquad  
{\cal L}^{-,}{}_l^a:=\rho_l^a(\R^{-1}{}^{(1)})\R^{-1}{}^{(2)} \label{frt}  
\ee  
by the relation $Z^h_k=(S{\cal L}^{-,}{}^h_i){\cal L}^{+,}{}^i_k$.
Eq.
(\ref{compact}) implies that $T$ is real, and
\be 
Z^h_k{}^{\star}=Z^k_h, \qquad\qquad 
(\c{L}^{\pm,}{}^i_j)^{\star}=S\c{L}^{\mp,}{}^j_i= 
g_{ih}\c{L}^{\mp,}{}^h_k g^{kj},             \label{qReal} 
\ee 
as $\rho$ is a $\star$-representation; the second  
equality in (\ref{qReal})$_2$ is based on the following   
useful property of the $N$-dimensional representation  
of $U_qso(N)$: 
\be  
\rho^a_b(Sh)=g^{ad}\rho^c_d(h)g_{cb}.             \label{Sonrho}  
\ee  

We recall that $U_qso(N)$ is a Ribbon Hopf algebra
\cite{ResTur90}:   the ribbon element $w\in U_qso(N)$
is a special, central element such that
\bea 
&& w^2=u_1 S(u_1), \qquad\qquad 
u_1:= (S\R^{(2)}) \R^{(1)}, \label{defw}\\ 
&& \Delta(w)=(w\otimes w)T^{-1} ,
\qquad\qquad Sw=w=S^{-1}w.           \label{coprodw} 
\eea 
It is well-known \cite{Dri90} that there exist isomorphisms
$U_hso(N)[[h]]\simeq Uso(N))[[h]]$ of $\star$-algebras over $\b{C}[[h]]$. This
essentially means that it is possible to express the elements of either
algebra as power series in $h=\ln q$ with coefficients in the other. In
particular $w$ has an extremely simple expression in terms of 
 the quadratic Casimir  $C$ of $so(N)$\footnote{This can be easily proved using 
the properties of the Drinfel'd twist $\F$ and the  
relation $\R=\F_{21}\, q^{X^a\!\otimes\! X_a}\, \F^{-1}$.}: 
\be
w=q^{-C}\!=\!e^{-hC}\!=\!\1\!+\!O(h), \qquad\qquad C\!:=\!X^aX_a=:L(L\!+\!N\!-\!2).
\ee
($\{X^a\}$ is a basis of $so(N)$).
We denote by $v:=w^{1/2}$, $\nu:=w^{1/4}$ and by $\tilde w,\tilde v,\tilde
\nu,\tilde T,\tilde Z$  the analogs of $w,\nu,T,Z$ obtained by replacing 
$\R$ by $\tilde\R$. As an immediate  consequence 
$$ 
\tilde w=q^{-C}q^{-\eta^2}, \qquad \tilde v=q^{-\frac
C2}q^{-\frac{\eta^2}2},\qquad\tilde \nu= q^{-\frac
C4}q^{-\frac{\eta^2}4},\qquad \tilde T= T q^{2\eta\otimes\eta}.  
$$ 
Since $C,-\eta^2$ are real (even positive-definite), if $q>0$ 
all these elements make sense either as positive-definite formal power series  
in $h$ of the form  $1+O(h)$, or as additional positive-definite generators 
of our Hopf $\star$-algebra. In Section \ref{PDHPO} we shall make them
into  positive-definite operators acting on the spaces of functions 
and of $p$-forms on $\b{R}_q^N$. 

All the information on the $\star$-algebras ${\cal DC}^*,H$ 
and the right action can be encoded in the 
cross-product $\star$-algebra ${\cal DC}^*\cocross H$.
We recall that this is $H\otimes {\cal DC}^*$
as a vector space,
and so we denote as usual $g\otimes a$ simply by $ga$;
that $H\1_{{\cal DC}^*}$, $\1_H{{\cal DC}^*}$ are subalgebras
isomorphic to $H, {{\cal DC}^*}$, and so we omit to 
write either unit $\1_{{\cal DC}^*},\1_H$
whenever multiplied by non-unit elements;  
that for any $a\in{\cal DC}^*$, $g\in H$ the product fulfills
\be
a g=g_{(1)}\, (a\tl g_{(2)}).                \label{crossprod}
\ee
${\cal DC}^* \cocross H$ is a $H$-module algebra  itself, if we
extend $\tl$ on $H$ as the adjoint action, namely as
$h\tl g= Sg_{(1)}\, h\, g_{(2)}$.
In view of (\ref{crossprod}),
this formula will correctly reproduce the action also on
the elements of ${\cal DC}^*$, and therefore on
{\it any} element $h\in {\cal DC}^*\cocross H$.
 The ``cross commutation relations'' (\ref{crossprod}) on the
generators $\sigma^h$ and $Z^i_j,\eta$ take the form
\bea
&&\sigma_1Z_2=\hat R_{12}Z_1\hat R_{12}\sigma_1
\qquad\qquad \mbox{i.e.}\qquad\qquad \sigma^hZ^i_j=
\hat R^{hi}_{km}Z^k_l\hat R^{lm}_{nj}\sigma^n,  \qquad\qquad 
                                         \label{crosscommrel}\\
&&x^i\eta=(\eta+1)x^i, \qquad\qquad
\xi^i\eta=(\eta+1)\xi^i, \qquad\qquad
\partial^i\eta=(\eta-1)\partial^i.\qquad\qquad 
\eea
The right relation in (\ref{crosscommrel}) is the translation of the left
one, where the conventional matrix tensor 
notation has been used.

\noindent
An alternative $\star$-structure for the 
whole ${\cal DC}^* \cocross H$ will be given in
(\ref{starsimid'}).

As shown in \cite{Fio95cmp,ChuZum95}, there exists a 
$\star$-algebra homomorphism 
\be 
\varphi: \A\cocross H\to \A,              \label{Hom1} 
\ee 
acting as the identity on $\A$ itself, 
\be 
\varphi(a)=a \qquad\qquad a\in \A,     \label{Hom1'} 
\ee 
where $H$ is the Hopf algebra $H=U_qso(N)$, and $\A={\cal H}$ 
is the deformed Heisenberg algebra. In \cite{Fio04} we have
extended $\varphi$ to the Hopf algebra 
$H=\widetilde{U_qso(N)}$ introducing an additional generator  
$\eta'=\varphi(\eta)\in{\cal DC}^*$ subject to the condition 
$\varphi(q^{\eta})=q^{\eta'}=q^{-N/2}\Lambda$, so that
\be 
[\eta',x^i]=-x^i \qquad 
[\eta',\partial^i]=\partial^i\qquad[\eta',\xi^i]=0 \qquad \eta'1|=q^{-N/2}.
\label{eta'drel}
\ee 
For real $q$, $\varphi$ is even a ${\star}$-algebra homomorphism.
Applying $\varphi$ to both sides of (\ref{crossprod}) one 
finds in particular 
\be 
a \,\varphi(g)=\varphi(g_{(1)})\, (a\tl g_{(2)}). \label{crossprod'} 
\ee 
In the sequel we shall often use the the shorthand notation 
\be 
\varphi(g)=:g',  \qquad \qquad \qquad g\in H.          \label{short} 
\ee 
We shall need in particular the images $Z'{}^h_k=\varphi(Z^h_k)$ 
explicitly. We  determine them here, starting from an Ansatz 
inspired by the images  $\varphi_l(Z^h_k)$ 
found in Ref. \cite{ChuZum95} for the analogous map 
$\varphi_l: U_qso(N)\cross {\cal H}\to {\cal H} $ 
(where $U_qso(N)$ acts with a {\it left} action): 
 
\begin{prop} Let $q\in\b{R}$.  
Under the $\star$-algebra map 
$\varphi: {\cal H}\cocross H\to {\cal H}$ the  
$\varphi(Z^h_k)$ are given by
\be 
Z'{}^h_k= 
q^{-2}\delta^h_k+ q^{-1}k\partial^hx^jg_{jk} 
-q^{-1-N}kx^h\hat\partial^jg_{jk}- \frac{k^2q^{-2}}{1+q^{N-2}} 
\partial^h r^2\hat\partial^jg_{jk}        \label{varphi} 
\ee 
where we have defined $k:=q-q^{-1}$. Moreover 
\be 
g'\1|=\epsilon(g)\1\qquad\qquad g\in H. 
\ee 
\label{prop1}
\end{prop} 
 The latter relation
together with (\ref{crossprod'}) implies  
\be 
g' f|=f\tl S^{-1}g.                     \label{imagedeco} 
\ee
In particular we find (\ref{vvalue}) on the spherical harmonics of level $l$.

One may ask if $\varphi$ trivially extends to a map 
of the type (\ref{Hom1}-\ref{Hom1'})  
with the Heisenberg algebra ${\cal H}$ replaced 
by the whole ${\cal DC}^*$. The answer is no: by using formula
(\ref{new})  one easily finds the commutation relation 
\be 
\xi_1\,Z'_2=\hat R^{-1}_{12}Z'_1\hat R_{12}\xi_1, 
\ee 
which differs from what one would obtain from  
(\ref{crosscommrel}) with $\sigma^i=\xi^i$ applying 
such a $\varphi$. Clearly, this formula 
holds also if we replace the matrix $Z'$ with any of its powers 
$Z'{}^h$. Now, note that
the $\xi^i$ commute with $\Lambda$, see (\ref{Lambdaprop})$_3$.
Recalling \cite{FadResTak89} that the center 
${\cal Z}(U_qso(N))$ of $U_qso(N)$ is
generated by the Casimirs $C_l$ defined by
\be 
C_l:=\mbox{tr}[UZ^h],\qquad\quad U^i_j:=g^{ik}g_{jk},\qquad\quad  
l=1,2,...,\left[N/2\right]  \qquad             \label{casimirs} 
\ee 
one easily checks and concludes that 
\be 
[\xi^i,C'_h]=0 \qquad \Rightarrow \qquad 
\left[\xi^i,\varphi\Big({\cal Z}(H)\Big)\right]=0 
\ee 
with $H=\widetilde{U_qso(N)}$, in particular $[\xi^i,\tilde w']=0$ , 
whereas $\xi^i$ do not  
commute with the center ${\cal Z}(H)$ itself 
(in fact $[\xi^i,C_h]\neq 0$, $[\xi^i,\eta]\neq 0$).

\subsection{Vielbein basis, Hodge map and Laplacian} 
\label{hodge}  

The set of $N$ exact forms $\{\xi^i\}$ is a natural basis 
for the  ${\cal H}$-bimodule 
${\cal DC}^1$, as well as for the  
 ${\cal H}\cocross \widetilde{U_qso(N)}$-bimodule 
${\cal DC}^1\cocross \widetilde{U_qso(N)}$.
In Ref.  \cite{CerFioMad00,Fio04} we
introduced  ``frame''  \cite{DimMad96} (or ``vielbein'') bases
$\{\theta^i \}$ and $\{\vartheta^i\}$ for the two, which are very useful for many purposes. 
These 1-forms are given by 
\bea 
&&\vartheta^i:=q^{-\eta-\frac N2}{\cal L}^{-,}{}^i_l\xi^l=  
\xi^m q^{1-\eta} \rho_m^j(u_4) {\cal L}^{-,}{}^i_j, 
                                         \label{defvarthetaa}\\ 
&&\theta^i:=\Lambda^{-1}\varphi({\cal L}^{-,}{}^i_l)\xi^l 
=\Lambda^{-1}\xi^h U^{-1}{}^i_k\varphi({\cal L}^{-,}{}^k_jU_h^j) 
= \Lambda^{-1}\xi^h \varphi(S^2{\cal L}^{-,}{}^i_h) \qquad\quad
\label{defthetaa'}   
\eea 
[$u_4:=\R^{-1}{}^{(1)}S^{-1}\R^{-1}{}^{(2)}$, and 
${\cal L}^{\pm,}{}^i_l$
are the FRT generators, see (\ref{frt})], and are characterized by the 
property
\be 
[\vartheta^i,{\cal H}]=0\qquad\qquad[\theta^i,{\cal H}]=0.  
                                              \label{framecond'} 
\ee 
They satisfy the same commutation relations as the $\xi^i$.
As already recalled,  from (\ref{xixirel}) it follows \cite{Fio94}
that $\mbox{dim}(\bigwedge^N)= 1$. The  
 matrix elements of the $q$-epsilon tensor are defined 
\cite{Fio94} up to a normalization constant $\gamma_N$ by
either relation
\be
\xi^{i_1}\xi^{i_2}...\xi^{i_N}= d^N\!x\:\varepsilon^{i_1i_2...i_N}, 
\qquad \qquad\theta^{i_1}\theta^{i_2}...\theta^{i_N}=  
dV\:\varepsilon^{i_1i_2...i_N},             \label{defqepsilon}
\ee
where
\be 
\gamma_N\,d^N\!x:=\xi^{-n}\xi^{1-n}...\xi^n\in\mbox{$\bigwedge^N$}, \qquad
\qquad \gamma_N\,dV:=\theta^{-n}\theta^{1-n}...\theta^n.    \label{defdNx} 
\ee 
One finds \cite{Fio04} that
the ``volume form'' $dV$ is  central in ${\cal DC}^*$ and equal to 
$dV=d^N\!x\Lambda^{-N}$. 
As a consequence of (\ref{Lambdaprop}), $dV|=d^N\!x$. 

Note that (\ref{framecond'}) in particular implies
$[\theta^i, \varphi(g)]=0$ for any $g\in U_qso(N)$. Going
to the differential basis $\xi^h$ by means of the
inverse transformation of (\ref{defthetaa'}) one
finds the following commutation relations between the $\xi^h$
and $g'=\varphi(g)$:
\be
\xi^h \varphi(g)= \varphi(S{\cal L}^{-,}{}^h_i \,g\,{\cal L}^{-,}{}^i_l)\xi^l
\qquad\quad
\varphi(g)\xi^h = \xi^l\varphi(S^2\!{\cal L}^{-,}{}^i_l\, g\, S\!{\cal
L}^{-,}{}^h_i).              \label{new}
\ee

As shown in \cite{Fio04}, for any $p=0,1,...,N$ one can define a
$U_qso(N)$-covariant,  ${\cal H}$-bilinear map 
\be 
*:{\cal DC}^p\to{\cal DC}^{N-p} 
\ee 
(the ``Hodge map''), such that ${}^*\1=dV$ and on each ${\cal DC}^p$  (and
therefore on the whole ${\cal DC}^*$)\footnote{
There is no sign at the rhs of (\ref{involution})  
[contrary to the standard $(-1)^{p(N-p)}$
of the undeformed case] because of the non-standard 
ordering of the indices in (\ref{defHodge1}). The latter in 
turn is the only correct one: had we used a different 
order, at the rhs of (\ref{involution}) tensor products 
of the matrices $U^{\pm 1}$, instead of the unit matrix, would have 
appeared, because of the property \cite{Ste96}  
\be  
\epsilon^{i_1 ... i_N} = (-1)^{N-1} U^{i_i}_{j_1}   
\epsilon^{i_2 ... i_N j_1}.                  \label{cycl_eps}  
\ee    
}
\be 
*^2\equiv *\circ *=\id       \label{involution} 
\ee 
by setting on the monomials in the $\theta^a$ 
\be 
{}^*(\theta^{a_1}\theta^{a_2}...\theta^{a_p}) 
=\,c_p\,\theta^{a_{p+1}}...\theta^{a_N} 
\varepsilon_{a_N...a_{p+1}}{}^{a_1...a_p},  \label{defHodge1} 
\ee 
(the normalization constants $c_p$ are given in \cite{Fio04}).
${\cal H}$-bilinearity of the Hodge map implies in particular 
\be 
{}^*(a\,\omega_p\, b)=a\,{}^*\omega_p\, b\qquad\qquad 
\forall \,a,b \in {\cal H}, \quad \omega_p\in{\cal DC}^p;  \label{Hbil} 
\ee 
i.e. applying Hodge and multiplying by ``functions 
or differential operators'' 
are commuting operations, in other words 
a differential form $\omega_p$ and its Hodge 
image have the same commutation 
relations with $x^i,\partial^j$. Restricting the domain of $*$ to the unital 
subalgebra $\widetilde{\Omega}^*\subset {\cal DC}^*$ generated by 
$x^i,\xi^j,\Lambda^{\pm 1}$ one obtains also a $U_qso(N)$-covariant,
$\widetilde{F}$-bilinear map 
\be
*:\widetilde{\Omega}^p\to\widetilde{\Omega}^{N-p}    \label{restri}
\ee
fulfilling again ${}^*\1=dV$  and (\ref{involution})  
(here $\widetilde{F}\equiv \widetilde{\Omega}^0$).  
The restriction (\ref{restri}) is the notion closest to the conventional notion
of a Hodge map on $\b{R}_q^N$: as a matter of fact, there is no
$F$-bilinear restriction of $*$ to $\Omega^*$. Note  
however that  $\widetilde{\Omega}^*$ is not closed 
under the $\star$-structure $\star$\footnote{In  Ref. \cite{CerFioMad00}  
we introduced a different $\star$-structure under which 
$\widetilde{\Omega}^*$   is closed.}. 

One would think that, since the vielbein $\theta^a$ do not belong
to $\Omega^*$, they cannot be used to describe a $p$-form
$\omega\in\Omega^*$ through components  
$\omega^{\theta}_{a_p...a_1}\in F$.
On the contrary, in section \ref{hilbert} we shall give a 
very useful notion of such components.

Finally, introducing the exterior coderivative
\be
\delta:= -{}^*\, d\,{}^*
\ee 
one finds that
on all of ${\cal DC}^*$, and in particular on all of $\Omega^*$,
the Laplacian $\Delta^{[\tilde\nu'{}^{-1}]}:=d\,\delta+\delta\,d$ is given by
\be
\Delta^{[\tilde\nu'{}^{-1}]}:=d\,\delta+\delta\,d=-q^2\,\partial\cdot\partial\Lambda^2=
-q^{-N}\hat\partial\cdot\hat\partial    
\label{Laplacian}
\ee
For the exterior coderivative $\hat\delta:=- {}^*\,\hat d\,{}^*$
of the ``hatted'' differential calculus
 one similarly finds
that the Laplacian 
$\Delta\equiv\Delta^{[\tilde\nu']}:=\hat d\,\hat \delta+\hat \delta\,\hat d$
is equal to $\Delta=-q^{-2}\,\hat\partial\cdot\hat\partial\Lambda^{-2}
=-q^N\partial\cdot\partial$.
The reason for the awkward superscripts ${}^{[\tilde\nu']}, {}^{[\tilde\nu'{}^{-1}]}$
will appear clear in section \ref{hilbert}.

\subsection{Integration over $\b{R}_q^N$ and naive scalar products} 
\label{Integration} 

In defining integration over  $\b{R}_q^N$, i.e. 
a suitable $\b{C}$-linear functional 
$$ 
f\in \Gamma\subset F\to  
\left(\int_q f\,d^N\!x\right) \,\in \b{C}, 
$$ 
we adopt the approach of 
Ref. \cite{Ste96} (already sketched in \cite{HebWei92}),  rather than the
preceding one of Ref.'s  \cite{Fio93,KemMaj94,fiothesis}\footnote{The
construction of  \cite{Fio93,KemMaj94,fiothesis} is purely algebraic, namely 
based on the fact that by repeated application of the Stokes theorem 
one can reduce $\int_q d^N\!x\,f$ to  
 $\int_q d^N\!x\,e_{q^2}[-r^2/a^2]$ for any function 
$f=e_{q^2}[-r^2/a^2]p(x)$ 
where $e_{q^2}[-r^2]$ is the $q$-gaussian and $p$ is 
a monomial in $x^i$; by linearity this can be extended also 
to power series $p(x)$ in a certain (not so large) class with fast 
decrease at infinity.}, since the former is applicable  
to a larger domain $\Gamma \subset F$ of ``functions'' 
(specified below). Going to ``polar coordinates''  
$\{x^i\}\to\{t^i,r\}$, $f(x)=f(t,r)$, allows to define the integral 
decomposing it into an integral over the ``angular coordinates'' $t^i$, i.e. 
over the $q$-sphere $S_q^{N\!-\!1}$, 
followed by the integral over the ``radial coordinate'' $r$: 
$$ 
\int_qf(x)\,d^N\!x=\int\limits^{\infty}_0dr\, m(r)\, 
r^{N\!-\!1} \int_{S_q^{N\!-\!1}}d^{N\!-\!1}\!t\,f(t,r). 
$$ 
Up to a normalization factor $A_N(q)$ 
(playing the role of the volume of $S_q^{N\!-\!1} $), 
which we here choose to be 1  
for the sake of brevity, 
the integration $\int_{S_q^{N\!-\!1}}d^{N\!-\!1}\!t$ coincides 
with the projection $f\in\Gamma\to f_0\in\Gamma_0$, 
where  
$\Gamma_0=\Gamma\cap \b{C}[[r,r^{-1}]]$ is the ``zero angular  
momentum'' subspace  of $\Gamma$ [see (\ref{ldeco})]: 
$\int_{S_q^{N\!-\!1}}d^{N\!-\!1}\!t\,f(x)=  f_0(r)$. This implies 
\be 
\int_qf(x)\,d^N\!x=\int\limits^{\infty}_0dr\,   m(r)\,
r^{N\!-\!1} f_0(r). 
                     \label{defint}\ee 
This has to be understood  as an integral of  
the {\it analytic continuation} of $f_0(r)$ to 
$\b{R}^+$, if $f_0$ is not assigned as a function on
$\b{R}^+$ from the very beginning; by $dr$ we mean 
Lebesgue measure, whereas $dr\, m(r)\equiv d\mu(r)$
denotes a Borel measure 
fulfilling the $q$-scaling property $ d\mu(qr)=q\,d\mu(r)$
(in other words the ``weight'' $m(r)$ fulfills $ m(qr)= m(r)$), 
which ensures the invariance under $q$-dilatations 
\be 
 \int_q f(qx) d^N\!(qx)\equiv\int_q \Lambda^{-1}f(x)| d^N\!(qx) 
= \int_q f(x) d^N\!x               \label{scaling_x}. 
\ee 
The ``weight'' $m_{J,r_0}(r):=|q-1|\sum_{n\in\b{Z}}r\delta(r-r_0q^n)$  
gives the  
socalled Jackson integral, $m(r)=1$ the standard Lebesgue integral, 
over $\b{R}^+$. 
Thus we can define integration on the functional space
$$
\Gamma=F\,\setminus\,\left\{f_0\in\b{C}[[r,r^{-1}]]\:\: |\:\:
\int_qf_0\,d^N\!x=\pm\infty \right\}. 
$$
For real $q$ integration over $\b{R}_q^N$ fulfills the following 
properties: 
\bea 
&& \left(\int_q f\,d^N\!x\right)^{\star} 
\:  = \int_q f^{\star}d^N\!x\qquad\qquad\qquad 
\qquad\qquad\qquad\quad\mbox{reality} \qquad 
 \label{reality_x} \\ &&\int_q 
f^{\star}f d^N\!x\quad \geq 0,  \qquad 
 \mbox{and =0 iff }f=0 
\qquad\qquad\qquad\mbox{positivity} \qquad  \label{pos_def_x} \\ 
&&\int_q \left(f\, d^N\!x\right)\tl g=\epsilon(g)\int_q f\,d^N\!x 
\qquad\qquad\quad\: U_qso(N)\mbox{-invariance}  \qquad    \label{cov} 
\eea 
Moreover, if $f$ is a regular function decreasing faster than 
$1/r^{N-1}$ as  $r\to\infty$ the Stokes 
 theorem holds 
\be 
\int_q \partial_i f(x)|\,d^N\!x=0, 
\qquad\qquad \int_q \hat\partial_i f(x)|\,d^N\!x=0. \label{stokes} 
\ee 
Properties (\ref{cov}-\ref{stokes}) express invariance respectively 
under deformed `infinitesimal translations and rotations'. 
On the contrary, the cyclic property for the integral of a product of functions
is $q$-deformed \cite{Ste96} .

Integration of functions immediately leads to integration 
of $N$-forms $\omega_N$. Upon moving all the $\xi$'s to the right 
of the $x$'s and using (\ref{defqepsilon}) we can express  
$\omega_N$ in the form $\omega_N=fd^N\!x$, and just have to set  
\be 
\int_q \omega_N=\int_q f\,d^N\!x. 
\ee 
Then eq. (\ref{stokes}) takes the form $\int_q d \omega_{N\!-\!1}|=0$, 
$\int_q \hat d \omega_{N\!-\!1}|=0$.
Finally, using Stokes theorem it is easy to show that for any $p=0,1,...N$ 
and any $\alpha_p\in{\cal DC}^p$, $\beta_{N\!-\!p}\in{\cal DC}^{N\!-\!p}$
\be
\int_q \alpha_p\, {}^{\star}\: \beta_{N-p}|=\int_q
(\alpha_p|)^{\star}\,\beta_{N-p}|                         \label{Leib} 
\ee
provided the product $\alpha_p\, {}^{\star}\beta_{N-p}|$ 
decreases fast enough as  $r\to\infty$. 
Because of the $\b{C}$-linearity of $\int_qd^N\!x$ and properties 
(\ref{reality_x}), (\ref{pos_def_x}), (\ref{Leib})  
one can introduce the (naive) scalar products of  two ``wave-functions''  
$\bphi,\bpsi\in F$ and more generally of two ``wave-forms'' 
$\balpha_p, \bbeta_p\in \Omega^p$ by  
\be
\la \bphi ,\bpsi\ra := \int_q \bphi^{\star} 
\bpsi \,d^N\!x,   \qquad      \qquad                        
\la\balpha_p, \bbeta_p\ra:=\int_q \balpha_p^{\star}\:{}^*\bbeta_p |.
                                        \label{basicscalprod}
\ee
{}From the decomposition (\ref{ldeco2}) for $\bphi,\bpsi$ and  
the orthonormality relations   
$\int_{S_q}d^{N\!-\!1}\!t\,S_l^I{}^{\star}S_{l'}^{I'}=  
(S_l^I{}^{\star}S_{l'}^{I'})_0=\delta_{ll'}\delta^{II'}$  
we find 
\bea  
\la \bphi ,\bpsi\ra &=&\int\limits^{\infty}_0dr\,  
r^{N\!-\!1} m(r)(\bphi^{\star}\bpsi)_0(r)\nn   
&=&\sum\limits_{l,l'=0}^{\infty}\sum\limits_{I,I'}  
(S_l^I{}^{\star}S_{l'}^{I'})_0\int\limits^{\infty}_0dr\,  r^{N\!-\!1}  
m(r)\phi_{l,I}^{\star}(r)\psi_{l',I'}(r)\nn  
&=&\sum\limits_{l=0}^{\infty}\sum\limits_I  
 \la \phi_{l,I},\psi_{l,I}\ra',         \label{interm} 
\eea
where we have introduced the `reduced scalar product'  
\be  
\la \phi,\psi\ra':=\int\limits^{\infty}_0dr\,   
r^{N\!-\!1} m(r)\phi^{\star}(r)\psi(r)= 
\int\limits^{\infty}_{-\infty}dy\,   
e^{Ny} \tilde m(y)\tilde \phi^{\star}(y)\tilde \psi(y)                  
\label{reduced}   
\ee  
of two functions $\phi(r),\psi(r)$ defined on the positive real  
line, and we have defined $y:=\log r$, $\tilde m(y):=m(e^y)$. 
A glance to (\ref{vvalue}) is sufficient to verify that
for any real $a$ the operator $w'{}^a$ 
(in particular $\nu'{}^{\pm 1}$)
is Hermitean w.r.t $\la~,~\ra$.

Using (\ref{starxid}), Stokes theorem (\ref{stokes}) and the analog of
(\ref{Leibniz}) for the $\hat \partial$-derivatives,  we find
that the  $p^{\alpha}$ are not Hermitean w.r.t.  $\la~,~\ra$, but
\cite{Fio95JMP}:
\be
\la \bphi ,p^{\alpha}\bpsi\ra=\int_q
\bphi^{\star} p^{\alpha}\bpsi |\,d^N\!x=
 \int_q(\hat p^{\alpha}\bphi|)^{\star}\bpsi\, d^N\!x=
\la\hat p^{\alpha} \bphi, \bpsi\ra ,        \label{hermi0}
\ee
with $\hat p^{\alpha} =-i\hat\partial^{\alpha}$. Using Stokes theorem,
in the appendix we show that (\ref{basicscalprod})$_2$ equals
\be
\ba{lll}
\la\balpha_p, \bbeta_p\ra  
& =& \frac{1}{c_{N\!-\!p}}\int_q \balpha^{\theta\,a_p...a_1}{}^{\star} 
\bbeta^{\theta\,a_p...a_1}{}\,d^N\!x  \\[8pt]
& =& \frac{1}{c_{N\!-\!p}}\int_q \balpha'{}^{\theta\,a_p...a_1}{}^{\star} 
\bbeta'{}^{\theta\,a_p...a_1}|{}\,d^N\!x  
\ea          \label{scalprodp}
\ee 
where  we have introduced the notation
\be 
\omega_p=\xi^{i_1}...\xi^{i_p}\omega_{i_p...i_1}(x) 
=\theta^{a_1}...\theta^{a_p}\omega'{}^{\theta}_{a_p...a_1}
=:\theta^{a_1}...\theta^{a_p}\omega_{a_p...a_1}^{\theta}(x)|
\ee 
for any $p$-form $\omega_p\in\Omega^p$. We shall call the functions
$\omega_{i_p...i_1},\omega_{a_p...a_1}^{\theta}$ (note: also
the latter belong to $F$, {\it not} to
${\cal H}$!) the components 
of the $p$-form $\omega_p\in\Omega^p$ respectively  
in the bases $\{\xi^i\},\{\theta^a\}$.
The $\omega_{a_p...a_1}^{\theta}$ {\it must not} be confused with
the components $\omega'{}^{\theta}_{a_p...a_1}$ of 
$\omega_p$ in the basis $\{\theta^a\}$,  defined above by 
$\omega_p=:\theta^{a_1}...\theta^{a_p}\omega'{}^{\theta}_{a_p...a_1}$
(without the final vertical bar);
the latter belong to ${\cal H}$, because 
$\theta^a\in{\cal DC}^*\setminus\Omega^*$!
Clearly $\omega_{a_p...a_1}^{\theta}=\omega'{}^{\theta}_{a_p...a_1}|$.

The above ``open-minded'' definition implies the 
following generalized notion of transformation of the components of a  
given differential $p$-form under the change of basis of 1-forms
$\xi^i\leftrightarrow \theta^a$:
\be
\ba{l}
\omega_{i_p...i_1}(x)= \Lambda^{-p}\,
\varphi\Big(S^2({\cal L}^{-,}{}_{i_p}^{a_p} 
...{\cal L}^{-,}{}_{i_1}^{a_1})\Big)\omega_{a_p...a_1}^{\theta}(x)|\\[6pt]
\omega_{a_p...a_1}^{\theta}(x)= \Lambda^p\,
\varphi\Big(S({\cal L}^{-,}{}_{a_p}^{i_p} 
...{\cal L}^{-,}{}_{a_1}^{i_1})\Big)\omega_{i_p...i_1}(x)| .
\ea
\ee

In the appendix we also show 
\be 
\la\balpha_p, \bbeta_p\ra =\la{}^*\balpha_p,  
{}^*\bbeta_p\ra.                                \label{bla2} 
\ee

Formula (\ref{scalprodp}) shows that (\ref{basicscalprod})$_2$ defines a 
``good'' scalar product in $\Omega^p$, reducing it
to the scalar product in $\bigwedge^p F$. 
In particular if $p=0$ then $\balpha_0,\bbeta_0\in F$ 
and we recover the scalar product (\ref{basicscalprod})$_1$, 
because 
$$
\int_q \balpha_0^{\star} \,{}^*\bbeta_0| 
=\int_q \balpha_0^{\star} \,dV\,\bbeta_0\,| 
=\int_q \balpha_0^{\star} \,\bbeta_0\,dV| 
=\int_q \balpha_0^{\star} \,\bbeta_0 \,d^N\!x .
$$
One defines a `naive'  Hilbert space of square integrable functions 
on $\b{R}_q^N$ by 
\be
\tilde L_2^m:=\left\{ {\bf f} (x)\equiv
\sum\limits_{l=0}^{\infty}\sum\limits_I S_l^I\,f_{l,I}(r)
\,\in F \:\: |\:\: \la {\bf f} ,{\bf f} \ra \, <\infty\right\}
\label{defL_2naive}
\ee 
(the superscript $m$ refers to the choice of
the radial measure $m$), and similarly one defines `naive' 
Hilbert space of square integrable  $p$-forms.

\sect{The $\star$-structure expressed by similarity
transformations}
\label{simi}

\begin{theorem} For positive $q$ 
the $\star$-structure of ${\cal DC}^*$
given in  (\ref{starx}-\ref{starxid}) can be
expressed in the form
\bea
x^i{}^{\star} &=& x^hg_{hi},\hfill \\
\xi^i{}^{\star}
&=& q^N\,\xi^hg_{hj}\,Z'{}^j_i\Lambda^{-2},
\label{starsimixi}\\
\partial^i{}^{\star}&=&-q^{\frac {1-N}2}  v'{}^{-1}\partial^hg_{hi}
v'\Lambda=- \tilde v'{}^{-1}\partial^h
\tilde v'{}g_{hi},   \label{starsimid}  \\                    
d^{\star} &=& - \tilde v'{}\,d\,\tilde v'{}^{-1}, \label{starsim}
\\         \theta^{\star} &=& \tilde w'{}\,\theta\,\tilde w'{}^{-1}.
\label{starsimith}  \eea
\label{starsimi}
\end{theorem}
\noindent
(The proof of the theorem is in the appendix.)
By the linear transformation $V^{\alpha}_i$ 
(see subsection \ref{preli1})we obtain a set of derivatives $\partial^{\alpha}$ such that on 
$-i\partial^{\alpha}$ 
 $\star$ acts as a similarity  
transformation:  
\be 
p^{\alpha}\equiv -i\partial^{\alpha}:=-iV^{\alpha}_i\partial^i 
\qquad\qquad \Rightarrow\qquad\qquad p^{\alpha}{}^{\star}= 
\tilde v'{}^{-1}p^{\alpha}\tilde v'{} 
\label{starsimil} 
\ee  


Incidentally, one can endow {\it the whole}
${\cal A} \cocross H$ with an alternative $\star$-structure
by keeping Eq. (\ref{starx})
unchanged  while removing
the map $\varphi$ from (\ref{starsimixi}-\ref{starsimid})
and readjusting the normalization factors in the latter formulae:
\be
\ba{lll}
\xi^i{}^{\star'}
&=& q^N\,\xi^hg_{hj}\,Z^j_iq^{-2\eta}
=\tilde w^{-1}\,\xi^hg_{hi}\,\tilde w, \label{starsimixi'}\\
\partial^i{}^{\star'}&=&- \tilde v^{-1}\partial^h
\tilde vg_{hi}, 
\\
d{}^{\star'} &=& -\tilde w^{-1}\,\xi^i\tilde v\partial_i \tilde v
=-\xi^i\tilde v\partial_i \tilde v^{-1}
\\         \theta^{\star'} &=& \tilde w\,\theta\,\tilde w^{-1};
\ea                       \label{starsimid'}
\ee
the second equality in the first line is easily proved
by means of the formulae given in Section \ref{vf}
and (\ref{wsigmarel}). We see that $\star'$ acts as a
similarity transformation also on the differentials 
$\xi^{\alpha}=V^{\alpha}_i\xi^i$.

Using (\ref{starsimid}), (\ref{starsimid'}), (\ref{frt}), (\ref{qReal}), the fact
that for real $q$ $\varphi$ is a $\star$-algebra map and
the relation $Z^h_k=(S{\cal L}^{-,}{}^h_i){\cal L}^{+,}{}^i_k$
  it is now straightforward to prove 

\begin{prop}

For real $q$ 
\be 
\vartheta^i{}^{\star'}=\vartheta^jg_{ji},\qquad\quad  
\theta^i{}^{\star}=\theta^jg_{ji},\qquad\quad
dV^{\star'}=dV=dV^{\star}.          \label{thetastar} 
\ee 
Moreover, the $\star$-structure and
the Hodge  map commute: 
\be 
({}^*\omega_p)^{\star}={}^*(\omega_p^{\star}). 
      \label{Hodgestar}  
\ee 
\end{prop}

\sect{New solutions for old problems:
improved real momentum, scalar products and Hermitean conjugation}
\label{hilbert}

We come now to some problems addressed in the introduction.

\begin{enumerate}
\item {\bf Quantum mechanics on $\b{R}_q^N$ as a configuration space.} One
question already  asked in the literature \cite{Wei94,fiothesis,Fio95JMP}
is: what is the ``right'' momentum sector subalgebra ${\cal P}$ 
within algebra of observables ${\cal H}$? In particular, what 
should be considered the ``right'' square momentum (i.e. Laplacian)
\cite{HebWei92,Wei94,fiothesis,Fio95JMP}? What are their spectral decompositions?

\item {\bf Field theory on $\b{R}_q^N$.} What is 
the ``right'' kinetic term in the action functional of a field-theoretic
model on $\b{R}_q^N$? This is clearly related also
to the question: what is the ``right'' propagator after quantization
of the model?
 
\end{enumerate}

As for problem 1., we wish to fulfill at least the
following  requirements. ${\cal P}$ must be: 
1. isomorphic to $F'$ (and therefore to $F$); 
2. closed under the action of $\widetilde{U_qso(N)}$; 
3. closed under the $\star$-structure. 
The solution proposed in \cite{Wei94,Fio95JMP} was essentially 
the subalgebra ${\cal P}\subset {\cal H}$ generated by 
the $p_R^{\alpha}$ defined by 
\be 
p_R^{2i+1}\!=\!\partial^i\!+\!\partial^i{}^{\star}\!=\! 
\partial^i\!-\!q^{-N}\hat\partial^jg_{ji}\qquad\quad 
p_R^{2i}\!=\!i[\partial^i\!-\!\partial^i{}^{\star}]= 
i[\partial^i\!+\!q^{-N}\hat\partial^jg_{ji}],     \label{realp'} 
\ee 
(where we adopt the indices' convention of \cite{Ogi92}, as in the previous
section)
 and in \cite{Fio95JMP} we even erroneously stated that it was  
uniquely determined 
(the proof of Theorem 2 of \cite{Fio95JMP} has a bug). 
The $p_R^{\alpha}$ are real and fulfill relations (\ref{ddrel}),
whereas (\ref{dxrel}), (\ref{dxirel})  are replaced by rather complicated
ones involving the angular momentum components [see relation (3) in 
\cite{Wei94} for the $\b{R}_q^3$ case]. Finding eigenfunctions of
a complete set of commuting observables including one or more
$p_R^{\alpha}$ is thus a rather hard task. Trying the same
even with just the square momentum (i.e. Laplacian)
$p_R^{\alpha}\cdot p_R^{\alpha}$ leads to lengthy calculations and complicated
formulae.\footnote{To see this, note that $p_R^{\alpha}\cdot p_R^{\alpha}$
is a combination of $\partial\cdot\partial$, $\hat\partial\cdot\hat\partial$ and
$\hat\partial\cdot\partial$. The latter in its own is an alternative, simpler
candidate for a real Laplacian, and in fact was diagonalized in Ref.
\cite{HebWei92}, formula (40),  where  a rather long
expression for its eigenvalues (involving
also the orbital angular momentum number $l$) was found. This is related to the
occurrence of the angular momentum in the commutation relations between
these Laplacians and the coordinates $x^i$.}

{}On the basis of the results of the previous section 
one could propose as an alternative solution 
that ${\cal P}\subset {\cal H}$ be the subalgebra generated by 
the  $\tilde p^{\alpha}$ defined by 
\be 
\tilde p^{\alpha}:=-iV^{\alpha}_i\tilde\partial^i,\qquad\qquad\qquad 
\tilde \partial^i:= \tilde \nu'{}^{-1}\partial^i 
\tilde \nu'.                            \label{realpt}
\ee 
Also the $\tilde p^{\alpha}$ are real. 
They fulfill relations (\ref{ddrel}), (\ref{dxirel}), 
whereas (\ref{dxrel}) is to be replaced by a so complicated one
that probably it cannot be put in closed form.\footnote{At least, 
one advantage is however that  the Laplacian  
$-\tilde p\cdot \tilde p\equiv\tilde\partial\cdot\tilde\partial$ 
is equal to $\partial\cdot\partial\Lambda q^{1-\frac N2}$ and therefore its
commutation relation with the coordinate $x^i$ is pretty manageable 
for iterated applications,
$$
\tilde\partial\cdot\tilde\partial\, x^i=(1\!+\!q^{2\!-\!N})
q^{-\frac N2}\partial^i\Lambda+qx^i\tilde\partial\cdot\tilde\partial,
$$
whereas the commutation relation of $-p_R\cdot p_R$ with $x^i$ 
is more complicated.} Similarly, one can introduce a purely 
imaginary nilpotent exterior derivative by 
\be 
\tilde d:=\tilde \nu'\, d\,\tilde \nu'{}^{-1} \qquad\qquad\Rightarrow 
\qquad\qquad \tilde d^{\star}=-\tilde d; 
\ee 
unpleasently it doesn't fulfill the ordinary Leibniz rule any more. 

As we now point out, the choice among the set 
$\{p_R^{\alpha}\}$, the  $\{\tilde p^{\alpha}\}$,   
the  $\{p^{\alpha}\}$ or any other  
set of derivatives, or between $\tilde d$ and $d$, will have physical 
significance only together with a specific choice of the scalar product  
within the Hilbert space upon they are meant to act. 
The standard `naive' scalar product (\ref{basicscalprod}) is just
{\it one} of the possible choices, but {\it not the only} one;
our goal is to adapt this choice to the choice of the (most
manageable) momentum components and exterior derivative.
Both the $p_R^{\alpha}$ and the $\tilde p^{\alpha}$ 
are (formally) Hermitean  w.r.t. the `naive' scalar product $\la~,~\ra$: 
\be
\la \bphi ,p_R^{\alpha}\bpsi\ra=\la p_R^{\alpha}\bphi , 
\bpsi\ra, \qquad\qquad
\la \bphi ,\tilde p^{\alpha}\bpsi\ra=\la \tilde p^{\alpha}\bphi , 
\bpsi\ra.                                  \label{ppreal}
\ee
The first equality (on the appropriate domains) follows from (\ref{hermi0}), 
and was already proved in \cite{Wei94,Fio95JMP,Wes99}; as we shall see
in section \ref{PDHPO},
the second actually holds (on the appropriate domains) if the radial
measure $m(r)$ is 1 or satisfies some other specific condition.
As already noted, the computation of
the action of either $p_R^{\alpha}$ or $\tilde p^{\alpha}$ is rather complicated
because none of them fulfills a simple
Leibniz rule like (\ref{Leibniz}).
As an alternative, we tentatively introduce  the `improved' scalar products
\be
(\cbphi ,\cbpsi) := \la \tilde \nu'{}^{-1}\cbphi, \tilde \nu'{}^{-1}\cbpsi\ra\qquad\qquad           
(\calpha_p, \cbeta_p):= \la  \tilde \nu'{}^{-1}\calpha_p,  \tilde \nu'{}^{-1}\cbeta_p\ra,
                                                       \label{scalprod}
\ee
the `improved'  Hilbert space of square integrable functions 
on $\b{R}_q^N$
\be
\check L_2^m:=\left\{{\bf f}(x)\equiv
\sum\limits_{l=0}^{\infty}\sum\limits_I S_l^I\,f_{l,I}(r)
\,\in F \:\: |\:\: ( {\bf f},{\bf f}) <\infty \right\},
\label{defL_2im}
\ee 
and similarly the `improved'  Hilbert space
of square integrable $p$-forms.
Under the conditions specified in Section \ref{PDHPO} 
the (in the algebraic sense) positive-definite  elements
$\tilde \nu'{}^{\pm 1}$ can be represented as Hermitean,
positive-definite  pseudodifferential operators on 
appropriate domains. Then 
\be
(\cbphi ,\cbpsi)=
\int_q \cbphi^{\star} \tilde \nu'{}^{-2} \cbpsi|\,d^N\!x, \qquad  \qquad                       
       (\calpha_p, \cbeta_p)=\int_q \calpha_p^{\star}\:\:{}^*\:\tilde
\nu'{}^{-2}\cbeta_p|.             
\ee
As a consequence of Theorem \ref{starsimi} and of the equality 
$ \tilde \nu'{}^2 =\tilde v' $ we obtain 
\be
(\calpha_p, \hat d\cbeta_{p\!-\!1}) = (\hat \delta\calpha_p,
\cbeta_{p\!-\!1}),  
 \qquad\qquad
(\hat d\cbeta_{p\!-\!1},\calpha_p) = (\cbeta_{p\!-\!1},\hat \delta\calpha_p)
\ee
and the (formal) hermiticity of both the
momenta $p^{\alpha}=i\partial^{\alpha}$ and the Laplacian
$\Delta$ w.r.t. the `improved' scalar product $(~,~)$: 
\be
(\cbphi ,p^{\alpha}\cbpsi) =   (p^{\alpha} \cbphi ,\cbpsi),   \qquad \qquad    
(\calpha_p,\Delta \cbeta_p) =(\Delta \calpha_p,\cbeta_p). \label{hermi!} 
\ee
In other words, the hermiticity of  $\tilde
p^{\alpha},\tilde\partial\cdot\tilde\partial$ 
w.r.t. 
$\la~,~\ra$ becomes equivalent to the hermiticity of
$p^{\alpha},\Delta\propto\partial\cdot\partial$  w.r.t. $(~,~)$!
If we impose the relation $\bphi = \tilde \nu'{}^{-1}\cbphi$ 
we can regard $\bphi$, $\cbphi$  as wave-functions  
representing the same ket and $\tilde p^{\alpha}$, $p^{\alpha}$ as  
pseudodifferential operators representing the same abstract operator 
in two different, but physically equivalent (configuration-space) `pictures',
because
\be
(\cbphi, \cbpsi)= \la\bphi, \bpsi\ra.
\ee 

Our answer to problem 1. is thefore as follows: In the original, `naive'
picture  the momentum
observables act on a wave-function $\bphi(x)$ 
as the pseudodifferential operators $\tilde p^{\alpha}$,
whereas the `position' observables 
act  simply by (left) multiplication by $x^{\alpha}$, yielding
$x^{\alpha}\bphi(x)$.  This picture is thus more convenient to compute the
action of the latter than the action of the former.  Instead in the second,
`improved'  picture the momentum operators act on a wave-function $\cbphi(x)$
as the differential operators $p^{\alpha}$, whereas the `position' observables,
act  as the pseudodifferential operators 
$\tilde \nu'x^{\alpha}\tilde\nu'{}^{-1}$.  
Therefore the second picture   is definetely more convenient for computing the
action of the momentum operators, as well as for answering
questions 2 (as we shall see below).

This notion of `picture' can be generalized as follows.   For any
pseudodifferential operator $\sigma=\id+O(h)$ depending only 
on $C',\eta'$, 
we introduce the ``$\sigma$-picture'' by 
\be 
\ba{l} 
{\bf f}^{[\sigma]}:= \sigma {\bf f}| \\[6pt] 
\la {\bf f},{\bf g}\ra^{[\sigma]} := \la \sigma^{-1}{\bf f},\sigma^{-1} {\bf
g}\ra   \equiv \int_q (\sigma^{-1}{\bf f}|)^{\star}\: 
\sigma^{-1}{\bf g}|\,d^N\!x \qquad \qquad \\ [6pt] 
{\cal O}^{[\sigma]}:=  \sigma \, {\cal O}\, \sigma^{-1} 
\ea 
\ee 
for ${\bf f},{\bf g}\in F$, ${\cal O}\in {\cal H}$ 
(note that for $\sigma=1$ one recovers the original  
picture). For our purposes it will be enough to stick 
to pseudodifferential operators of the form $\sigma
=q^{a(\eta'\!+\!b){}^2}g(C)$, where $b$ is a real constant and $g(C)$ is a
positive-definite pseudodifferential operator depending only on the quadratic
Casimir of $so(N)$. We tentatively introduce the
``Hilbert space  of square integrable functions  
on $\b{R}_q^N$ in the $\sigma$-picture'' by  
\be 
\tilde L_2^{m,\sigma}:=\left\{{\bf f}(x)\equiv 
\sum\limits_{l=0}^{\infty}\sum\limits_I S_l^I\,f_{l,I}(r) 
\,\in F \:\: |\:\: \Vert {\bf f}\Vert^2_{\sigma}<\infty\right\}, 
\label{defL_2s} 
\ee  
where $\Vert {\bf f}\Vert^2_{\sigma}:=\la {\bf f},{\bf f}\ra^{[\sigma]}$.  In particular, 
$\cbphi=\bphi^{[\tilde\nu']}$, $\bphi=\bphi^{[1]}$,
$\check L_2^m=\tilde L_2^{m,\tilde\nu'}$.  
Then, trivially 
\be 
\ba{l} 
\bphi^{[\sigma]}\!\in\! \tilde L_2^{m,\sigma}
\quad\Leftrightarrow \quad\bphi\!\in \!\tilde L_2^m \\[8pt] 
\la \bphi^{[\sigma]} ,\bpsi^{[\sigma]}\ra^{[\sigma]}=\la \bphi ,\bpsi\ra, 
\ea 
\ee 
and [denoting by $D^{[\sigma]}\big({\cal O}^{[\sigma]}\big)$ 
the domain of operator ${\cal O}^{[\sigma]}$ within $\tilde L_2^{m,\sigma}$] 
\be  
\ba{l} 
\bphi^{[\sigma]}\in  D^{[\sigma]}\Big({\cal O}^{[\sigma]}\Big)\subset\tilde L_2^{m,\sigma} 
\quad\Leftrightarrow \quad 
\bphi\in D({\cal O})\subset\tilde L_2^m,\\[8pt] 
{\cal O}^{[\sigma]}\bphi^{[\sigma]}|=({\cal O}\bphi|)^{[\sigma]}, 
\ea 
\ee  
implying that one can describe the same ``physics''  
by any of the $\sigma$-pictures. So one can choose the most convenient
for each computation.

\noindent
The generalization of the notion of $\sigma$-pictures to forms is straightforward. 

In section \ref{PDHPO} we determine radial measures $m$ and 
for each  $\sigma$ of the above type a ($m$-dependent) 
subspace $L_2^{m,\sigma}\subset \tilde L_2^{m,\sigma}$  
and define $\sigma$ as a pseudodifferential 
operator such that 
 \be 
\la {\bf f} , {\bf g}\ra^{[\sigma]}= \la {\bf f} ,{\bf
g}^{[(\sigma\sigma^{\star})^{-1}]}\ra  =\la {\bf f}^{[(\sigma\sigma^{\star})^{-1}]} ,
{\bf g}\ra            \label{hermi2}  
\ee 
for any ${\bf f},{\bf g}\in L_2^{m,\sigma}$, in particular 
\be 
(\cbphi,\cbpsi) = \la\hbphi,\cbpsi\ra = \la\cbphi,\hbpsi\ra,
                                               \label{hermis} 
\ee 
where $\hbphi\equiv\bphi^{[\tilde\nu'{}^{-1}]}$,
$\cbphi\equiv \bphi^{[\tilde\nu']}$.
After the replacements
$\bphi\to\hbphi$, $\bpsi\to\cbpsi$, (\ref{hermi0}) becomes 
\be 
\la \hbphi , p^{\alpha}\cbpsi\ra=  
\la \hat p^{\alpha} \hbphi ,\cbpsi\ra.    \label{hermi1} 
\ee 
Then (\ref{starsimid}), (\ref{realpt}) will imply 
(\ref{ppreal})$_2$, (\ref{hermi!})$_1$ 
respectively for any $\bphi, \bpsi\in D(\tilde p^{\alpha})$, 
and $\cbphi, \cbpsi\in D^{[\tilde\nu']}(p^{\alpha})$ 
(note that with our notation  
$p^{\alpha}=\tilde p^{\alpha}{}^{[\tilde\nu']}$) and more generally 
\be 
\la \bphi^{[\sigma]}  \!,\!\tilde p^{\alpha}{}^{[\sigma]}\bpsi^{[\sigma]} \ra^{[\sigma]} = 
\la\tilde p^{\alpha}{}^{[\sigma]}\bphi^{[\sigma]}  ,\bpsi^{[\sigma]} \ra^{[\sigma]} 
\ee 
for any $\sigma$ and
$\bphi^{[\sigma]}, \bpsi^{[\sigma]}\in D^{[\sigma]}(\tilde p^{\alpha}{}^{[\sigma]})$. 
 

 
\medskip
 
As an application, we recall how one can diagonalize observables 
of ${\cal P}$ using improved pictures. 
In Ref. \cite{Fio95JMP} we constructed irreducible 
$\star$-representations of  the $\star$-algebra 
${\cal P}\cocross H\subset {\cal H}$ and diagonalized within the latter a
complete set of commuting observables,
consisting not only of the square total momentum
 $P\cdot P=:(P\cdot P)_n$,but of all the $(P\cdot P)_a:=\sum_{j=-a}^{a}P^jP_j$ with $a=1,2,...,n$
(these are the squares of the projections of the
 momentum on the hyperplaneswith coordinates $P^{-a},P^{1-a},...,P^a$), of $P^0$ (only
 for odd $N$), andof the generators $K^a$ of the Cartan subalgebra
 of $U_qso(N)$.Diagonalization 
 was performed first at the abstract level, i.e. eigenvectorswere abstract kets
 and ${\cal P}$ was the $\star$-algebra generated byabstract
 $U_qso(N)$-covariant generators $P_i$ fulfilling (\ref{ddrel})and the same $\star$-relations (\ref{starx}) as the $x_i$.
Then we realized the scheme in $\b{R}_q^N$-configuration space in two different 
realizations, i.e. pictures: 
in the first one (which we called ``unbarred'') $P_i$
were realized  as 
$-i\Lambda\partial_i=-i\tau\tilde\partial_i\tau^{-1}$,
in the second (which we called ``barred'') 
 the $P_i$ were realized as$-i\hat\partial_i\Lambda^{-1}=\tau^{\star}{}^{-1}\tilde\partial_i\tau^{\star}$
where $\tau:=\nu'q^{(\eta'\!+\!N\!+\!1)^2/4}$.
In the previous notation they amount respectively to the
$\sigma=\tau$ and the
$\sigma=\tau^{\star}{}^{-1}$ pictures\footnote{We warn the reader that
in the conventions of \cite{Fio95JMP} $\Lambda$ is what here is denoted by
$\Lambda^{-1}$, and conversely.}. To
compute the action
 of $P_i$ either one is much more convenient than the
`naive' one, where $P_i$ are realized as the pseudodifferential operators
$-i\tilde \partial_i$,
 because of the relatively simple commutation relations
(\ref{dxrel}), (\ref{Lambdaprop}) and the analogous ones involving the
$\hat\partial_i$.  For   $0<q<1$ we found the
following spectral decompositions of the above  observables: 
\be
\ba{l} 
(p\cdot p)\, \bphi^{[\tau]}_{\bpi ,{\bf j}} = \kappa^2 q^{2\pi_n}\,\bphi^{[\tau]}_{\bpi ,{\bf j}}, 
\\[8pt] 
(p\cdot p)_a\, \bphi^{[\tau]}_{\bpi ,{\bf j}} = \kappa_a^2
q^{\sum_{k=a}^n2\pi_k}\,\bphi^{[\tau]}_{\bpi ,{\bf j}}, \label{eigen}  \\[8pt]
 p_0\,  \bphi^{[\tau]}_{\bpi,{\bf j}} =\kappa_0 q^{\pi_0}\, \bphi^{[\tau]}_{\bpi ,{\bf j}} 
\qquad\qquad \mbox{ (only for odd $N$)}; 
\ea
\ee
here $\kappa\equiv \kappa_n$ is a positive constant chacterizing the 
irreducible representation
(by a redefinition of $\pi_n$ it can be always chosen in $[1,q[$), and
$$ 
\kappa_a=\kappa q^{n\!-\!a}
\sqrt{\frac{1\!+\!q^{-2\rho_a}}{1\!+\!q^{N\!-\!2}}}, 
\qquad \qquad \kappa_0=
\pm \kappa q^n \sqrt{\frac{1\!+\!q^{-1}}{1\!+\!q^{N\!-\!2}}} \:\mbox{ (only
for odd $N$)},  
$$
whereas $\bpi,{\bf j}$ are
vectors (the component $j_a$ of ${\bf j}$ labels eigenvalues
of $K^a$) with suitable  \cite{Fio95JMP}
integer components, in particular $\pi_n\in\b{Z}$ and $\pi_h\in\b{N}$ if $h<n$.    
Up to normalization,  in the unbarred realization (or 'picture') 
the eigenfuntions
$\bphi^{[\tau]}_{\bpi,{\bf j}}$ with $\bpi=
{\bf 0}$ will be given by \cite{Fio95JMP}
 $$
\bphi^{[\tau]}_{{\bf 0} ,{\bf j}} \sim (x^{-n})^{j_n}\!...(x^{-2})^{j_2}\cdot
\left\{\ba{lr} 
(x^{-1})^{j_1}e_{q^{-1}}[i\kappa_0x^0] \qquad &\mbox{if }
N\!=\!2n\!+\!1 \\[8pt]
 (x^{-sign(j_1)\cdot1})^{|j_1|}\varphi_{q^{-1}}^J\left(x^1 x_1\frac
{q\kappa{}^2}{q^{2\!-\!N}\!+\!1}\right)
 \quad &\mbox{if } N\!=\!2n
\ea\right.
 $$
where $J:=\sum_{a=1}^nj_a$ and, having set $(l)_q:=(q^l\!-\!1)/(q\!-\!1)$,
\be
e_q(z) :=\sum\limits_{l=0}^{\infty}{z^l \over (l)_q!}, 
\qquad\varphi_q^J(z):=\sum\limits_{l=0}^{\infty}\frac{(-z)^l}
{(l)_{q^2}!(l\!+\!J)_{q^2}!}.            \label{specf}
\ee
(As we expect, for odd $N$ in the limit $q=1$ $\bphi^{[\tau]}_{{\bf 0} ,{\bf
0}}$ formally becomes a plane wave orthogonal to the $x^0$ coordinate).  The $\bphi^{[\tau]}_{{\bf 0} ,{\bf j}}$ can be also obtained
from the cyclic eigenfunction $\bphi^{[\tau]}_{{\bf 0} ,{\bf 0}}$ by applying
to the latter suitable elements in ${\cal P}\cocross H$. The
$\bphi^{[\tau]}_{\bpi,{\bf j}}$ with $\bpi\neq {\bf 0} $ are obtained
applying to $\bphi^{[\tau]}_{{\bf 0},{\bf j}}$  powers of the
$\Lambda\partial_i$ with $i>0$.  We thus find relatively `tractable' 
eigenfunctions, which can be actually expressed through $q$-special functions
(see section \ref{defpseudodiff'}). Formula (\ref{eigen}) shows that these
operators have   very simple discrete spectra, essentially  consisting of
integer powers of $q$. As a matter of fact, the eigenfunctions are also
normalizable: this was proved in \cite{Fio95JMP} adopting a slightly different
definition of integration, and is true also adopting the definition of
integration \cite{Ste96}  recalled in section \ref{Integration}.\footnote{In
either case, the question of the normalizability of all
$\bphi^{[\tau]}_{\bpi,{\bf j}}$  is reduced to the question of the
normalizability of the cyclic eigenfunction $\bphi^{[\tau]}_{{\bf 0} ,{\bf
0}}$ by manipulations involving the use of Stokes theorem, similarly as in the
undeformed context the normalizability of the H\'ermite functions is reduced
to that of the  gaussian $e^{-r^2/2}$. That $\bphi^{[\tau]}_{{\bf 0} ,{\bf 0}}$
is normalizable is true by the definition of integration of \cite{Fio95JMP} in
the first case, and can be proved by a rather lenghty computation 
in the present case.}  
This situation is to be contrasted with the undeformed one,
where the corresponding operators have continuous spectra and generalized
eigenfunctions. Therefore $q$-deformation can be seen as a `regularizing'
device! Moreover, in section \ref{PDHPO}  we shall see that the constant
$\kappa{}$ characterizing the irreducible representation can take any value
if we choose a trivial radial weight [$m(r)\equiv 1$] in
(\ref{defint}), whereas (at least for even $N$) is {\it quantized} to a
specific value (defined up to powers of $q$)
if we choose a nontrivial $m(r)$. In other words,  in
the latter case the nature of space(time) fixes an energy scale independent of
the particular irreducible representation we have chosen, namely of the
particular type of particles we describe by the latter!

Similarly one can treat the case $q>1$. 

\bigskip
We come now to question 2.
The kinetic term in the action for a $p$-form (i.e. an antisymmetric
tensor with $p$-indices) Euclidean field theory with mass $M$ can be most simply introduced as
$$
{\cal S}_k=\Big( (\Delta+ M^2)  \calpha_k,\calpha_k\Big).
$$
It will be rather `tractable' because $\Delta=-q^N\partial\cdot\partial$ has
the rather 
 simple action (\ref{Boxxrel}) as a differential operator.Consider in particular
 a scalar field (i.e. $k=0$). The `propagator'(or Green function) $G(y,x)$ of the theory should be expressible
in terms of any orthonormal basis $\{\cbphi_{\pi_n,l,I} \}$ of eigenfunctions
of $\Delta+ M^2$,
 $\nu'$
\be
\ba{l}
(\Delta+ M^2) \cbphi_{\pi_n,l,I}= (\kappa^2q^{2\pi_n} +M^2) 
\cbphi_{\pi_n,l,I},
 \\[8pt]\nu'\cbphi_{\pi_n,l,I}=
q^{-l(l\!+\!N\!-\!2)/4} \cbphi_{\pi_n,l,I},
\ea                                                   \label{eigeneq}
\ee
and some other observables (whose eigenvalues we label by a multi-index $I$)
commuting with each other and making up a 
complete set, through the relatively simple formula
\bea
G(y,x)  &=& \sum\limits_{\pi_n,l,I} \cbphi_{\pi_n,l,I}(y) 
\left[\tilde\nu'{}^{-2}(\Delta+ M^2)^{-1}\cbphi_{\pi_n,l,I}|\right]^{\star}(x)\nn
&=& \sum\limits_{\pi_n,l,I} \cbphi_{\pi_n,l,I}(y)\frac
1{\kappa^2q^{2\pi_n}\!+\!M^2} 
 [\tilde\nu'{}^{-2}\cbphi_{\pi_n,l,I}|]^{\star}(x)\nn
&=&\sum\limits_{\pi_n,l,I}
\cbphi_{\pi_n,l,I}(y)\frac{q^{l(l\!+\!N\!-\!2)/2}}{\kappa^2q^{2\pi_n}\!+\!M^2} 
[q^{\eta'{}^2/2} \cbphi_{\pi_n,l,I} | ]^{\star}(x),
\eea
where $y^i$ denote the generators of another copy of $\b{R}_q^N$.
If we choose $I$ as the multi-index
 labelling spherical harmonics(\ref{defS}) one thus looks for the basis elements in the 
form  $\cbphi_{\pi_n,l,I}= S_l^I \phi_{\pi_n,l}(r)$. 
Using the formulae given in appendix \ref{related} reduces Eq.
(\ref{eigeneq})$_1$ to a $q$-difference equation for $\phi_{\pi_n,l}(r)$;
solving it is now an affordable task, which is left as a job for
future work.

\sect{Defining the pseudodifferential operators $q^{a(\eta'\!+\!b)^2}$}
\label{PDHPO}

As said, in order that the formal
considerations of the previous
section are implemented  at the operator level we  have to make sense
out of $\sigma =q^{a(\eta'\!+\!b)^2}g(C)$  as 
pseudodifferential operators on  
$F$ (more generally on
$\Omega^*$) and investigate whether we need to
restrict $\tilde L_{2,p}^{m,\sigma}$  to some subspace $L_{2,p}^{m,\sigma}$
in order that on the latter (\ref{hermi2}) holds. 
We are going to do this next, distinguishing the
case $m\equiv 1$ from the others. Clearly it is sufficient to do this
for $p=0$-forms, i.e. functions, because the form components are 
functions themselves.
Recalling the decomposition (\ref{ldeco2}) for $\bphi$, 
(\ref{vvalue}) 
and (\ref{interm}) we see that $g(C)$ fulfillsthe requirement, so the problem is 
reduced  to showing that one can define $q^{a(\eta'\!+\!b)^2}$ so that the latter
also does.
To define the action of 
$q^{a(\eta'\!+\!b)^2}$ on the functions $\phi_{l,I}(r)$
we perform the change of variable 
$r\to y:=\ln r$, whereby $\eta'= -\partial_y- N/2$
and $r^{N-1}dr=e^{Ny}dy$, for any function
$\phi(r)$ denote $\tilde \phi(y):=\phi(e^y)$, 
and express   $e^{yN/2}\tilde
\phi(y)=r^{N/2}\phi(r)$  in terms of its Fourier transform $\hat
\phi(\omega)$:  
\be 
e^{\frac N2 y}\tilde \phi(y)=\frac 1 
{\sqrt{2 \pi}}\int\limits^{\infty}_{-\infty} \hat \phi(\omega) e^{i\omega 
y}d\omega.                         \label{inversefourier}
\ee 
Here we are assuming in addition that  all
$e^{\frac N2 y}\tilde \phi_{l,I}(y)\in L_2(\b{R})\equiv L_2(\b{R},dy)$, 
in other words that
 all $\phi_{l,I}(r)\in L_2(\b{R}^+,dr^N)$, what guarantees that
the Fourier transform exists and is invertible.
One initial motivation behind such a change of variable is that $y$ is more
suitable to describe the behaviour of functions occurring in $q$-analysis,
notably $q$-special functions (which are tipically involved as
solutions of $q$-difference equations) as $r\to 0,\infty$ (i.e.
$y\to
 -\infty,\infty$), since often they wildly fluctuate 
as $r\to 0$ or as $r\to \infty$; this can be inferred from the typical
exponential scaling laws of 
the zeroes/poles $r_n$
 of $q$-special functions either 
as $r\to 0$ or $r\to \infty$\footnote{This happens for instance
with the $q$-gaussian $e_{q^2}[-r^2]:={}_0\varphi_0[q^2,(q^2-1)r^2]$:  
property (\ref{billo}) implies
$e_{q^2}[-q^2r^2]=[1-(q^2-1)r^2]e_{q^2}[-r^2]$, whence we see that  for 
$q>1$ and sufficiently large $r$ the modulus of
$e_{q^2}[-q^{2n}r^2]$ grows with $n$ and its sign flips at each step $n\to
n\!+\!1$.}. From (\ref{inversefourier})  we find
$$
e^{\frac N2 y}q^{a(\eta'\!+\!b)^2}\tilde \phi(y)|=q^{a\partial_y^2}e^{\frac N2
y}\tilde \phi(y)|=\frac 1  {\sqrt{2\pi}}\int\limits^{\infty}_{-\infty}d\omega
\hat \phi(\omega)   e^{i\omega y} q^{-a(\omega\!+\!ib)^2},
$$
i.e. $q^{a(\eta'\!+\!b)^2}$ acts as multiplication by $q^{-a\omega^2}$ 
on the Fourier transform, implying 
\be
\sigma\bphi(x)|=\frac {e^{-\frac N2 y}}{\sqrt{2\pi}}
\sum\limits_{l=0}^{\infty}g[l(l\!+\!N\!-\!2)]
\sum\limits_IS_l^I  \!\!\int
\limits^{\infty}_{-\infty}d\omega \hat \phi_{l,I}(\omega)  
 e^{i\omega y}q^{-a(\omega\!+\!ib)^2}.                \quad \label{b1}
\ee
Of course this is well-defined only
for $q^{-a}$ such that the integrals are.
We also easily see that one can extend the domain
of the partial derivatives $\partial^i$
to $\bphi$ with  $\phi_{l,I}(r)\in L_2(\b{R}^+,dr^N)$
using (\ref{aa1}) and (\ref{dir^2rel}),  provided we can extend also 
the action $\Lambda^{\pm 1}f(x)=f(q^{\mp 1}x)$ of $\Lambda^{\pm 1}\equiv
e^{\mp h\partial_y}$  on such $\bphi$'s; this is done of  course by setting 
\be
\Lambda^{\pm 1}\bphi(x)|=\frac {e^{-\frac N2 (y\mp h)}}{\sqrt{2\pi}}
\sum\limits_{l=0}^{\infty}
\sum\limits_IS_l^I  \!\!\int 
\limits^{\infty}_{-\infty}d\omega \hat \phi_{l,I}(\omega)   
 e^{i\omega (y\mp h)}.                \quad \label{bb1} 
\ee

In terms of Fourier transforms the reduced scalar product
(\ref{reduced})   becomes 
\be 
\la \phi ,\psi\ra'=  
\int\limits^{\infty}_{-\infty}d\omega'   
\int\limits^{\infty}_{-\infty}d\omega \:\hat \phi^{\star}(\omega)   
\hat\psi(\omega') \int\limits^{\infty}_{-\infty}\frac {dy}{2\pi}  
\tilde m(y)e^{i(\omega'-\omega) y}.                  \label{redscal'}
\ee 

\subsection{The case $m\equiv 1$}
\label{defpseudodiff}

In the case 
$m(r)=\tilde m(y)\equiv 1$ the third integral at the
rhs(\ref{redscal'}) reduces to   $\delta(\omega-\omega')$, implying 
\bea
\la \phi ,\psi\ra'&=&
\int\limits^{\infty}_{-\infty}d\omega\,   \hat \phi^{\star}(\omega)   
\hat\psi(\omega),\nn 
\la \bphi ,\bpsi\ra &=&\sum\limits_{l=0}^{\infty}\sum\limits_I 
\int\limits^{\infty}_{-\infty}d\omega\,   \hat \phi_{l,I}^{\star}(\omega)   
\hat\psi_{l,I}(\omega).
\eea
For $\bphi^{[\sigma]},\bpsi^{[\sigma]}\in F$,  this and (\ref{b1})
for $\sigma =q^{a(\eta'\!+\!b)^2}g(C)$ imply
\bea
&&\la\bphi^{[\sigma]} ,\bpsi^{[\sigma]}\ra^{[\sigma]}:=
\la \sigma^{-1}\bphi^{[\sigma]} , \sigma^{-1}\bpsi^{[\sigma]}\ra\nn
&& = \sum\limits_{l=0}^{\infty}\sum\limits_I g^2\big(l(l\!+\!N\!-\!2)\big)
 \int\limits^{\infty}_{-\infty}d\omega\,   
\Big(\hat\phi^{[\sigma]}_{l,I}(\omega)\Big)^{\star}   
\hat\psi^{[\sigma]}_{l,I}(\omega)q^{2ab^2-
2a\omega^2}\qquad\qquad\label{hermitic1}\\  
&&= \la \bphi^{[\sigma]} ,
\bpsi^{[\sigma^{\star}{}^{-1}]}\ra=\la \bphi^{[\sigma^{\star}{}^{-1}]} ,
\bpsi^{[\sigma]}\ra, \label{hermitic2}
\eea
in particular
\be
\Vert \bphi^{[\sigma]}\Vert^2_{\sigma}
=\sum\limits_{l=0}^{\infty}  g^2\big(l(l\!+\!N\!-\!2)\big)\sum\limits_I  
\int\limits^{\infty}_{-\infty}d\omega\,   |\hat \phi^{[\sigma]}_{l,I}(\omega)|^2 
q^{2ab^2-2a\omega^2} .      \qquad \quad      \label{domaincond}    
\ee
The function $\bphi^{[\sigma]}$ will belong to $\tilde L_2^{1,\sigma}$ if 
this is finite. If both $\bphi^{[\sigma]},\bpsi^{[\sigma]}\in \tilde
L_2^{1,\sigma}$ then by Schwarz inequality the rhs(\ref{hermitic1}) is finite as
well; then equalities in (\ref{hermitic2}) are just the proof of
relation (\ref{hermi2}) we were seeking for.

Note that in the present $m\equiv 1$ case by (\ref{interm}) 
the condition 
$\Vert\bphi^{[\sigma]}\Vert^2_{\sigma} <\infty$  characterizing
$\tilde L_2^{1,\sigma}$ implies  
$q^{a \eta'{}^2}\phi_{l,I}\in L_2(\b{R}^+,dr^N)$ for
all $l,I$, whence the assumed exixtence and invertibility of the
Fourier transform automatically follows.  
We summarize the results by stating the following

\begin{theorem}
If $m\equiv 1$, for any real $s$  the scalar product of the
Hilbert space $L_2^{1,\sigma}:=\tilde L_2^{1,\sigma}$ can be expressed by any
of the
 expressions in (\ref{hermi2}) and the $\tilde p^{\alpha}{}^{[s]}$ are 
(formally) hermitean
operators defined on $L_2^{1,\sigma}$.
\end{theorem}

\noindent
{\bf Remark.} If $q^a>1$ the factor
$q^{-2a\omega^2}$ in (\ref{domaincond})  acts as a ``UV regulator''.

\subsection{The case $m\neq 1$}
\label{defpseudodiff'}

The measure $m\equiv 1$ describes a continuous and
homogeneous space along the radial direction. It is important to 
leave room for a discretized space by allowing for a non-unit 
$m$, notably a measure concentrated in points, like Jackson's measure
$m_{J,r_0}(r)dr^N$, where
$$
m_{J,r_0}(r):=|q-1|\sum_{l\in\b{Z}}r\delta(r\!-\!r_0q^l)=|q-1|\sum_{l\in\b{Z}}
\delta(y\!-\!y_0\!-\!lh)
$$
(here $y_0=\log r_0$). The case $m\neq 1$ actually reveals to be rather 
interesting and rich of surprises; in the sequel we disclose
some of its features by performing a preliminary analysis,
leaving an exhaustive investigation as the subject for 
some other work.

We assume that all the $\phi_{l,I}(r)$ can be analytically 
continued to the complex $r$-plane. Sticking for simplicity
to the case that $\phi_{l,I}(r)$ are uni-valued,
the analytic
 continuation of $\tilde\phi_{l,I}(y):=\phi_{l,I}(e^y)$ 
will fulfill the periodicity
 condition
\be
\tilde \phi_{l,I}(y)=\tilde \phi_{l,I}\left(y+i2\pi\frac k{\gamma}\right),  
\qquad \qquad  k\in\b{Z}        \label{iperiodic}
\ee
with ${\gamma}=1$; more generally, they will also fulfill this condition
with ${\gamma}=2,3,...$ if $\phi_{l,I}(r)$ can be expressed in the form
$\phi_{l,I}(r)=\u{\phi}_{l,I}(r^{\gamma})$, with $\u{\phi}_{l,I}(z)$
uni-valued. Below we shall occasionally
suppress the subscripts $l,I$ in the intermediate
results to avoid a too heavy notation.
Now we compute the Fourier transform 
$\hat\phi$ of $\tilde \phi(y)e^{\frac N2  y}$ 
\be 
\hat \phi(\omega)=
\int\limits^{\infty}_{-\infty} \frac {dy}{\sqrt{2\pi}} 
\Phi(\omega,y),\qquad\qquad
\Phi(\omega,y):=\tilde \phi(y)e^{\frac N2 y-i\omega y}
                                         \label{Fourtransf}
\ee 
using the method of residues.
We first assume that $\phi(r)$ has no poles on $\b{R}^+$
(or equivalently that $\tilde \phi(y)$ has no poles on the 
real axis).
For $\omega<0$ the exponential $e^{-i\omega y}$ rapidly goes
to zero as $\Im(y)\to\infty$. Choose a contour  like the one
depicted in Fig 1, with $M\in\b{N}$. By (\ref{iperiodic}), 
the integral on the upper
horizontal side equals $-e^{iNM\pi+\omega M2\pi}$ times
(\ref{Fourtransf}), and therefore vanishes in the limit 
$M\to\infty$, together with the integral on the vertical sides.
Therefore, taking this limit
 we find
$$
\hat \phi(\omega)=i\sqrt{2\pi}\sum\limits_{\mbox{\tiny poles} 
\: y'\in \b{C}^+}
 \mbox{Res}\, \Phi(\omega,y')
=i\sqrt{2\pi}\sum\limits_{\mbox{\tiny poles} \: y'\in \b{C}^+} 
 \mbox{Res}\,\tilde\phi(y')e^{(\frac N2 -i\omega)y'}.
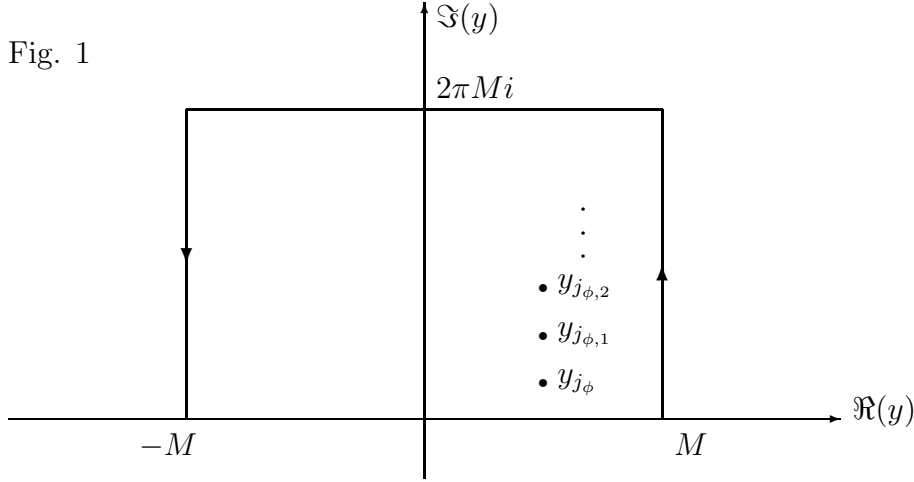
$$

\setlength{\unitlength}{2947sp}%
\begingroup\makeatletter\ifx\SetFigFont\undefined%
\gdef\SetFigFont#1#2#3#4#5{%
  \reset@font\fontsize{#1}{#2pt}%
  \fontfamily{#3}\fontseries{#4}\fontshape{#5}%
  \selectfont}%
\fi\endgroup%
\begin{picture}(7012,5000)(-4500,-4000) 
\thinlines 
\put(-4500,-3000){\vector( 1, 0){7000}} 
\put(-1000,-3500){\vector( 0, 1){4000}} 
\thicklines 
\put(-3000,-1700){\line( 0, -1){1300}} 
\put(1000,-1700){\line( 0, 1){1300}} 
\put(1000,-3000){\vector( 0, 1){1300}} 
\put(-3000,-400){\vector( 0, -1){1300}} 
\put(-3000,-400){\line( 1, 0){4000}} 
\put(1100,-3300){\makebox(0,0)[lb]{\smash{\SetFigFont{12}{14.4} 
{\rmdefault}{\mddefault}{\updefault}
$M$
}}}  
\put(-3400,-3300){\makebox(0,0)[lb]{\smash{\SetFigFont{12}{14.4} 
{\rmdefault}{\mddefault}{\updefault}
$-M$
}}}  
\put(-900,-300){\makebox(0,0)[lb]{\smash{\SetFigFont{12}{14.4} 
{\rmdefault}{\mddefault}{\updefault}
$2\pi Mi$
}}}  
\put(-900,300){\makebox(0,0)[lb]{\smash{\SetFigFont{12}{14.4} 
{\rmdefault}{\mddefault}{\updefault}
$\Im(y)$
}}}  
\put(2600,-3000){\makebox(0,0)[lb]{\smash{\SetFigFont{12}{14.4} 
{\rmdefault}{\mddefault}{\updefault}
$\Re(y)$
}}}  
\put(-4500,0){\makebox(0,0)[lb]{\smash{\SetFigFont{12}{14.4} 
{\rmdefault}{\mddefault}{\updefault}
$\mbox{Fig. 1}$
}}}  
\put(0,-2700){\circle*{60}}
\put(0,-2700){\makebox(0,0)[lb]{\smash{\SetFigFont{12}{14.4} 
{\rmdefault}{\mddefault}{\updefault}
$\:\: y_{j_{\phi}}$
}}}  
\put(0,-2300){\circle*{60}}
\put(0,-2300){\makebox(0,0)[lb]{\smash{\SetFigFont{12}{14.4} 
{\rmdefault}{\mddefault}{\updefault}
$\:\: y_{j_{\phi,1}}$
}}}  
\put(0,-1900){\circle*{60}}
\put(0,-1900){\makebox(0,0)[lb]{\smash{\SetFigFont{12}{14.4} 
{\rmdefault}{\mddefault}{\updefault}
$\:\: y_{j_{\phi,2}}$
}}}  
\put(300,-1700){\makebox(0,0)[lb]{\smash{\SetFigFont{12}{14.4} 
{\rmdefault}{\mddefault}{\updefault}
$\cdot$
}}}  
\put(300,-1500){\makebox(0,0)[lb]{\smash{\SetFigFont{12}{14.4} 
{\rmdefault}{\mddefault}{\updefault}
$\cdot$
}}}  
\put(300,-1300){\makebox(0,0)[lb]{\smash{\SetFigFont{12}{14.4} 
{\rmdefault}{\mddefault}{\updefault}
$\cdot$
}}}  
\end{picture}  

By (\ref{iperiodic}) the poles of $\tilde \phi(y)$,  and  
therefore of $\Phi(\omega,y)$, can be parametrized in the form  
\be 
y'_{j_{\phi},k}= y_{j_{\phi}}+2\pi \frac k{\gamma} i
\qquad\qquad 0<\Im(y_{j_{\phi}}) <2\frac{\pi}{\gamma} \label{poles1} 
\ee 
where $k\in\b{Z}$ and $j_{\phi}$ is some possible 
additional index. Therefore 
\bea 
\hat \phi(\omega)&=&i\sqrt{2\pi}\sum\limits_{j_{\phi}} 
\mbox{Res }\tilde \phi(y_{j_{\phi}})
e^{(\frac N2 -i\omega)y_{j_{\phi}}}\sum\limits_{k=0}^{\infty} 
e^{(\frac N2 -i\omega)i2\pi\frac k{\gamma}}\nn 
&=&\frac {i\sqrt{2\pi}}{1- e^{\frac{\pi}{\gamma}(iN+2\omega)}}                       
\sum\limits_{j_{\phi}}  \mbox{Res }\tilde \phi(y_{j_{\phi}})
e^{(\frac N2 -i\omega)y_{j_{\phi}}}           \label{explfour}
\eea 
since by (\ref{iperiodic}) 
$\mbox{Res} \,\tilde \phi(y_{j_{\phi}})
=\mbox{Res} \,\tilde \phi(y_{j_{\phi}}+i2\pi k/{\gamma})$. 
By applying the method of residues instead to an analogous  
clockwise contour 
in the lower complex $y$-half-plane $\b{C}^-$ one  
finds that the latter formula gives $\hat \phi(\omega)$ also for $\omega>0$. 

Note that if $N/{\gamma}$ is an even integer 
$\hat \phi(\omega)$ has a first order pole in $\omega=0$
 and
$\int^{\infty}_{-\infty}d\omega$ in (\ref{inversefourier}) 
has to be understood as a principal value integral 
around $\omega=0$, unless cancellations of contributions
of different poles $j_{\phi}$ occur.

Replacing (\ref{explfour}) in (\ref{b1}), 
if no $\phi_{l,I}(r)$ has poles on
$\b{R}^+$  we find
\be 
\sigma\bphi(x)|\!=\!i
\sum\limits_{l=0}^{\infty}g[l(l\!+\!N\!-\!2)] 
\sum\limits_IS_l^I  \sum\limits_{j_{l,I}} 
\mbox{Res }\tilde \phi(y_{j_{l,I}})\!\!\int 
\limits^{\infty}_{-\infty}d\omega\!
 \frac{e^{(i\omega -\frac N2)(y-y_{j_{l,I}})}
q^{-a(\omega\!+\!ib)^2}}{1- e^{\frac{\pi}{\gamma}(2\omega\!+\!iN)}}.    \quad
\label{b2}
\ee 
where we have used the short-hand notation 
$j_{l,I}:=j_{\tilde\phi_{l,I}}$. 
The integral is well-defined for $q^{-a}\le 1$,
i.e. $ah\ge 0$.
Note that if $ah>0$, because of the damping factor 
$q^{-a\omega^2}$, $\tilde w'{}^{-a}\bphi(x)$ 
has no more poles in $y=y_{j_{l,I}}$. Formula
(\ref{bb1}) will still give the action of 
$\Lambda^{\pm 1}$ on $\bphi$.

Let us now evaluate 
$\la\sigma\bphi,\sigma'\bpsi\ra$
(with $ah,a'h\ge 0$)
in the present case. By (\ref{interm}) and the previous
equation we find
\be
\la\sigma\bphi,\sigma'\bpsi\ra =
\sum\limits_{l=0}^{\infty}\sum\limits_I
gg'[l(l\!+\!N\!-\!2)] 
\la q^{a(\eta'\!+\!b)^2}\phi_{l,I},q^{a'(\eta'\!+\!b')^2}\psi_{l,I}\ra'
\ee
and 
\bea
&& \la q^{a(\eta'\!+\!b)^2}\phi,q^{a'(\eta'\!+\!b')^2}\psi\ra'
=\sum\limits_{j,j'} 
\left[\mbox{Res}\,\tilde \phi\left(y_{j}\right)\right]^{\star}
\mbox{Res}\,\tilde \psi\left(y_{j'}\right)\:M^j_{j'}(a\!,\!b;a'\!,\!b'),\label{redsc}\\
&&M^j_{j'}:=\int\limits^{\infty}_{-\infty}\!  d\omega'{} 
\int\limits^{\infty}_{-\infty}\!  d\omega\frac{e^{\frac
N2(y_{j}^{\star}\!+\!y_{j'})+ i(\omega
y^{\star}_{j}\!-\! \omega' y_{j'})} 
q^{-a(\omega\!+\!ib)^2\!-\!a'(\omega\!-\!ib')^2}}
{[1- e^{\frac{\pi}{\gamma}(2\omega\!-\!iN)}][1- e^{\frac{\pi}{\gamma}(2\omega'\!+\!iN)}]}
 \int\limits^{\infty}_{-\infty} \!
dy \,  \tilde m(y)e^{i(\omega'\!-\!\omega) y}\nonumber. 
\eea
[here $y_{j'}$  denote the pole locations of $\tilde\psi(y)$
with $0<\Im(y_{j'})<2\pi/{\gamma}$]. 
We ask whether 
$\la \phi, q^{a(\eta'\!+\!b)^2}\psi\ra'=\la q^{a(\eta'\!-\!b)^2}\phi,\psi\ra'$
for $\phi,\psi$ within a suitable space  of
functions to be identified. For $\gamma\in\b{N}$ and 
$\beta=0,\frac 12$  
 let 
\bea
&& \!\!L_2^{m,[\beta,\gamma]}\!:=\!
\Big\{\phi\!\in\! L_2\big(\b{R}^+\!,m(\!r\!)dr^N\!\big)\:\:\: |\:\:\:
\phi(r)=f(r)\u{\phi}(r^{\gamma}),\:\mbox{ where }\nn
&& \qquad\qquad \u{\phi}\mbox{ is
analytic with poles only in } z\!=\!-\!q^{n(j\!+\!\beta)}\!, 
\:j\!\in\!\b{Z} \Big\}. \qquad\quad\label{poles3}  
\eea
The poles of $\phi(r)$ will be only in
\be
r_{j,k}:=q^{j\!+\!\beta}e^{i\frac{\pi(2k\!+\!1)}{\gamma}} \label{poles4}
\ee
with $k=0,1,...,\gamma\!-\!1$ and $j$ belongs to some subset
$J\subset\b{Z}$,
and those of $\tilde\phi(y)$ only in
\be
y_{j,k}:=h(j\!+\!\beta)\!+\!i\frac{\pi(2k\!+\!1) }{\gamma}.
\ee
Condition (\ref{poles4}) amounts to saying that 
the pole locations lie on $\gamma$ special straight half-lines 
starting from $r=0$ and forming with each 
other angles equal to $2\pi/{\gamma}$, 
and are such that their absolute values are either
$q^j$ or  $q^{j+\frac 12}$,  with $j\in J\subset\b{Z}$.
The condition appearing in (\ref{poles3}) thus implies 
(\ref{poles1}) (with $\Im(y_j)=\pi/{\gamma}$),
whence (\ref{b2}-\ref{redsc}).
Thus, if  $\phi,\psi\,\in\, L_2^{m,[\beta,\gamma]}$, 
then 
\bea
M^j_{j'}(a\!,\!b;a'\!,\!b') &=&\frac 14\!\!\int\limits^{\infty}_{-\infty}\! 
d\omega'{}  \int\limits^{\infty}_{-\infty}\! 
d\omega \frac{q^{-a(\omega\!+\!ib)^2\!-\!a'(\omega\!-\!ib')^2\!+\!\frac N2
(j\!+\!j'\!+\!2\beta)}}{\sin[\frac{\pi}{2{\gamma}}(\!N\!-\!2i\omega)]
\sin[\frac{\pi}{2{\gamma}}(\!N\!+\!2i\omega')]} \times \nn 
&&\qquad\int\limits^{\infty}_{-\infty} \! dy \, \tilde m(y)
e^{i[(\omega'\!-\!\omega)(y\!-\!h\beta)\!+\!\omega jh\!-\!\omega'j'h]}.
\label{I3}   
\eea
Note that in (\ref{redsc}) one can consider the indices $j,j'$ as
running over the whole $\b{Z}$ for any $\phi,\psi$ because
the residues will vanish in the $y_j$ which are not poles for these
functions. Then one can consider $M(a\!,\!b;a'\!,\!b')$ as an universal infinite matrix
and express the lhs(\ref{redsc}) in terms of the row-by-column
matrix product
\be
\la q^{a(\eta'\!+\!b)^2}\phi,q^{a'(\eta'\!+\!b')^2}\psi\ra'
=R_{\phi}^{\dagger}\: M(a\!,\!b;a'\!,\!b')\: R_{\psi},             \label{phipsi'}
\ee
where by $R_{\phi}$ we have denoted the column vector with
infinitely many components 
$R_{\phi}^j$, $j\in\b{Z}$,
given by 
$R_{\phi}^j=\mbox{Res}\,\tilde\phi\vert_{y=\left[h(j\!+\!\beta)\!+\!i\pi/{\gamma}\right]}$.

Now, performing the change of integration variables  
$\omega\to -\omega'$ one immediately finds that
$M^{j'}_j(a'\!,\!b';a\!,\!b)=M^j_{j'}(a\!,\!b;a'\!,\!b')$. Moreover, taking
the complex conjugate and
 performing the change of integration variables  
$\omega\to -\omega$, $\omega'\to -\omega'$ we find that the
$M^j_{j'}$ are real,
\be
\Big[M^j_{j'}(a\!,\!b;a'\!,\!b')\Big]^{\star}=M^j_{j'}(a\!,\!b;a'\!,\!b').
\ee
By the $q$-scaling  property,    
the transformed weight $\tilde m(y):= m(e^y)$   is periodic  
with period $h=\ln q$; we shall also assume that  
$m$ is invariant under 
$r$-inversion\footnote{For the Jackson weight $m_{J,r_0}$ given 
above this  necessarily requires $r_0=1$ or $r_0=q^{1/2}$.}, 
so for any $k\in\b{Z}$  
\be 
\ba{lll}
&m(q^kr)=m(r) ,\qquad\qquad &m(r^{-1})=m(r),\\
\mbox{i.e. }\:&\tilde m(y\!+\!kh)=\tilde m(y),\qquad\qquad &\tilde m(-y)=\tilde
m(y). \ea                                       \label{qscaleinv} 
\ee 
Performing the change of integration variables 
$\omega'{}\leftrightarrow\omega$, 
$y\to -y+(j\!+\!j'\!+\!2\beta+2a'b'-2ab)h$ we now find
\be
M^j_{j'}(a\!,\!b;a'\!,\!b') = M^j_{j'}(a'\!,\!b';a\!,\!b), \qquad\mbox{ if
 }\,N/{\gamma}\!\in\!\b{N},\quad 2(a'b'\!-\!ab)\!\in\!\b{Z};       \label{exch} 
\ee
in fact, the  weight $\tilde m$ and the
last integral in (\ref{I3}) are automatically invariant under this
change of integration variables, whereas the condition $N/{\gamma}\in\b{N}$ ensures that
also the denominator in the first two is. From these relations we find
that the matrix $M$ is Hermitean:
\be
M^{\dagger}(a\!,\!b;a'\!,\!b')=M(a\!,\!b;a'\!,\!b').                  \label{Herm}
\ee
This is true in particular if $a=a'$, $b=b'$.
Choosing instead $a'=0=b'$  relations (\ref{exch}) and (\ref{Herm}) 
together with (\ref{redsc}) respectively imply
\bea
&&\la \phi, q^{a(\eta'\!+\!b)^2}\psi\ra'=\la q^{a(\eta'\!-\!b)^2}\phi,
\psi\ra',\label{hermi} \\ 
&&\la \phi, q^{a(\eta'\!-\!b)^2}\psi\ra'{}^{\star}=
\la q^{a(\eta'\!-\!b)^2}\psi,\phi\ra'=\la \psi, q^{a(\eta'\!+\!b)^2}\phi\ra'.
\label{sesqui} 
\eea
In formula (\ref{posicon}) in the appendix we give a necessary and sufficient
condition on the weight $m$ (which is satisfied in particular
by the Jackson measure) and on the parameters $a,h,\gamma$
in order that the
 positivity condition 
\be
\la \phi, q^{a\eta'{}^2}\phi\ra'\ge 0, \qquad\qquad    
\la \phi, q^{a\eta'{}^2}\phi\ra'=0 \quad\mbox{iff }\phi=0     \label{positi}
\ee
is fulfilled. We need this to be true
with any $a$ such that $ah\ge 0$, in particular with $a=1/2$ for
(\ref{hermi!})$_1$ to be valid,  or alternatively with
$a=-1/2$ for the analog of (\ref{hermi!})$_1$ with $p^{\alpha}$
replaced by the $\hat p^{\alpha}$ to be valid.
Then for any $\sigma=g(C)q^{a(\eta'\!+\!b)^2}$ with $2ab\in\b{Z}$ 
\be
\ba{rcl}
\la\bphi,\bpsi\ra^{[\sigma]} &=&
\sum\limits_{l,I} g^2[l(l\!+\!N\!-\!2)]
\la \phi_{l,I},q^{2a (\eta'{}^2\!+\!b^2)}\psi_{l,I}\ra'=
\la\bphi,\sigma^{-2}\bpsi\ra\\[8pt]
&=&\sum\limits_{l,I} g^2[l(l\!+\!N\!-\!2)]
\la q^{2a (\eta'{}^2\!+\!b^2)}\phi_{l,I},\psi_{l,I}\ra'
 =\la\sigma^{-2}\bphi,\bpsi\ra
\ea                                               \label{scalprod2}
\ee
defines a ``good'' scalar product within the 
the following subspace of $\tilde L_2^{m,s}$,
\be
L_2^{m,\sigma,[\beta,\gamma]}:=  
\{\bphi\equiv\sum\limits_{l,I}
 S_l^I\phi_{l,I}(r)\:\,
|\:\, \phi_{l,I}\in L_2^{m,[\beta,\gamma]}\:\,
\mbox{with}
\:\,\Vert \bphi\Vert_{\sigma}<\infty \} \label{defL_2_beta}
\ee
(here $\Vert \bphi\Vert_{\sigma}:=\la\bphi,\bphi\ra^{[\sigma]}$), 
making the latter a pre-Hilbert space. Relation (\ref{sesqui}) ensures
the sesquilinearity of $\la\:,\:\ra^{[\sigma]}$, (\ref{positi}) its positivity.
The $\tilde p^{\alpha}{}^{[\sigma]}$ are (formally) hermitean
operators on their domain within $L_2^{m,\sigma,[\beta,\gamma]}$, as a
consequence of (\ref{hermi}).  Investigating their essential self-adjointness
in the completed Hilbert space
is left as a job for future work. We collect the results by
stating the following 

\begin{theorem}
 Let $\beta\in\{0,1/2\}$, $\gamma\in\b{N}$ be a submultiple of
$N$, $ah\ge 0$, $4ab\in\b{Z}$, $\sigma=g(C)q^{a(\eta'\!+\!b)^2}$. 
Assume that the radial weight
$m(r)$ fulfills
 (\ref{qscaleinv}) and (\ref{posicon}), where 
$\check m(y)\equiv m\left(e^{h(y/2\!+\!\beta)}\right)$.
Then  (\ref{scalprod2}) defines the scalar product of a
pre-Hilbert space $L_2^{m,\sigma,[\beta,\gamma]}$  and the $\tilde
p^{\alpha}{}^{[\sigma]}$ are  (formally) hermitean
operators defined on $L_2^{m,\sigma,[\beta,\gamma]}$.
\end{theorem}

The spaces introduced in (\ref{defL_2_beta}) are very
interesting. Functions $\phi_{l,I}$ fulfilling (\ref{poles4}) are for
instance
\be
\frac 1{1+(q^{j\!+\!\beta}r)^{\gamma}}, \qquad \qquad f(r)\prod\limits_{l}\frac
1{1+(q^{j_l\!+\!\beta}r)^{\gamma}},                   \label{nice}
\ee
where $j_l\in\b{Z}$, $\beta=0,1/2$
and $f(r)$ is a polynomial or more generally analytic in a domain
including all $\b{R}^+$.
To this category belong also some $q$-special functions
 with distinguished (i.e. quantized) values of the
parameters characterizing them. Essentially
all special functions can be defined as particular cases of the
$q$-hypergeometric functions
${}_r\varphi_s(a_1,...,a_r;b_1,...,b_s;q,z)$\footnote{See for instance
\cite{GasRah90,KliSch97}, defined as (analytic continuations in the complex
$z$-plane of) 
\bea
&& {}_r\varphi_s(a_1,...,a_r;b_1,...,b_s;q,z):=\\
&& \sum\limits_{n=0}^{\infty}\frac{(a_1;q)_n...(a_r;q)_n}
{(b_1;q)_n...(b_s;q)_n}\left( (-1)^n 
q^{n(n\!-\!1)/2}\right)^{1\!+\!s\!-\!r} \frac{z^n}{(q;q)_n}\nonumber
\eea
(with parameters such that the series has at least a finite convergence
radius), where 
$$
(a;q)_0:=1,  \qquad\qquad
(a;q)_n:=\prod\limits_{i=0}^{n\!-\!1}(1-aq^i),\qquad n=1,2,...
$$
(whenever $|q|<1$ the latter definition makes sense also for $n=\infty$).

For instance the functions introduced in (\ref{specf}) can be expressed as
$$
e_q(z)={}_0\varphi_0(q,(1\!-\!q)z),\qquad
\varphi^J_q(z)=\frac 1{(J)_{q^2}!} \:
{}_2\varphi_1\left(0,0;q^{2J};q^2,-(1\!-\!q^2)^2z\right). 
$$ 
One can rewrite them in the form (\ref{nice})$_2$, using their interesting
properties (see e.g. \cite{KliSch97}}). For example
\bea
&&{}_0\varphi_0(q,z)=\prod\limits_{i=0}^{\infty}\frac 1{1-zq^i}
\qquad\qquad {}_1\varphi_0(a;q,z)=\prod\limits_{i=0}^{\infty}\frac
{1-a zq^i}{1-zq^i}             \qquad       \label{billo}
\\&&{}_2\varphi_1(a_1,a_2;b;q,z)=\frac{(a_2;q)_{\infty}(a_1z;q)_{\infty}}
{(b;q)_{\infty}(z;q)_{\infty}}\:
{}_2\varphi_1\left(\frac b{a_2},z;a_1z;q,a_2\right)    \qquad\label{billo'} 
\eea

Using (\ref{billo})$_1$ and (\ref{billo'})with $a_1,a_2=0$, $b=q^l$
($l\in\b{Z}$) one can check\footnote{Details will be given elsewhere} that the
eigenfunctions $\bphi^{[\tau]}_{\bpi ,{\bf j}}$ written in section
\ref{hilbert} belong to the space  $L_2^{m,\sigma,[\beta,\gamma]}$ where
$\sigma=\tau:=\nu'q^{(\eta'\!+\!N\!+\!1)^2/4}$, $\beta=0$, $\gamma=1$ provided
the energy scale $\kappa^2$ appearing in their definition is {\it quantized}
(up to powers of $q$) as follows: 
\be
\kappa^2=\frac{(1\!+\!q^{2\!-\!N})^2}{(1\!-\!q^2)^2}. 
\ee

\app{Appendix}

\subsection{Proof of Theorem \ref{starsimi} and related lemmas}
\label{related}

For   $\sigma^i=x^i,\xi^i,\partial^i$ we easily find
\be 
\sigma^i\tl w^{\pm 1}= q^{\pm(1\!-\!N)}\sigma^i, \qquad\qquad 
\sigma^i\tl \tilde w^{\pm 1}= q^{\mp N}\sigma^i.   \label{wonx} 
\ee 
In fact
\bea
&& \sigma^i\tl u_1 \stackrel{(\ref{fundrep})}{=}\sigma^j\rho^i_j(u_1 )
\stackrel{(\ref{defw})}{=}
\sigma^j\rho^i_h(S\R^{(2)}) \rho^h_j(\R^{(1)})\nn
&&\stackrel{(\ref{Sonrho})}{=} 
\sigma^j g^{il}\rho^m_l(\R^{(2)})g_{mh}\rho^h_j(\R^{(1)})
 =  \sigma^j g^{il}g_{mh}\hat R^{mh}_{jl}
\stackrel{(\ref{projectorR}),(\ref{Pt})}{=}
q^{1-N}\sigma^j g^{il}g_{jl}.
 \nonumber
\eea
Recalling (\ref{Sonrho}) 
and (\ref{defw}) we find $\sigma^i \tl w^2=\rho^i_j(u_1Su_1)=q^{2-2N} \sigma^i$
whence the first part of the claim.
The proof of the second statement is completely analogous. 
It is not difficult to check that (\ref{wonx}) implies 
\be 
w\,\sigma^i\,w^{-1}=q^{N\!-\!1}Z^i_j \sigma^j.       \label{wsigmarel}
\ee

\begin{lemma} Let $w_l:=q^{-l(l+N-2)}$. Then on the spherical harmonics of level $l$
(with $l=0,1,2,...$) 
\be 
w'S_l^I|=w_l\,S_l^I=S_l^I\tl w,  
\qquad\qquad w'{}^aS_l^I|=(w_l)^a\,S_l^I=S_l^I\tl w^a
\label{vvalue} 
\ee 
for any real $a$. In particular $ \nu'S_l^I| =q^{-l(l+N-2)/4}\,S_l^I=S_l^I\tl
\nu$ 
\end{lemma} 
\bp{} 
We determine the eigenvalue $w_l$ applying the pseudodifferential 
operator $w'\equiv \varphi(w)$ to $S_l^{n...n}=(t^n)^l$: 
\bea 
w' (t^n)^l| 
&\stackrel{(\ref{imagedeco})}{=}&(t^n)^l\tl S^{-1}w  
\stackrel{(\ref{coprodw})_2}{=} (t^n)^l\tl w\nn  
 &\stackrel{(\ref{modalg2})}{=}&  
t^n\tl w_{(1)}[(t^n)^{l-1}\tl w_{(2)}] \nn 
&\stackrel{(\ref{coprodw})_1}{=}& 
t^n\tl w T^{-1(1)}[(t^n)^{l-1}\tl wT^{-1(2)}] \nn 
&\stackrel{(\ref{wonx})_1}{=}& q^{1-N}w_{l-1}
t^n\tl T^{-1(1)}[(t^n)^{l-1}\tl T^{-1(2)}] \nn
&\stackrel{(\ref{modalg2})}{=}&q^{1-N}w_{l-1} 
(t^n\tl T^{-1(1)}) (t^n\tl T^{-1(2)}_{(1)})  ...
(t^n\tl T^{-1(2)}_{(l-1)})\nn 
&\stackrel{(\ref{fundrep})}{=}& q^{1-N}w_{l-1} 
\rho^n_{i_1} (T^{-1(1)}) \rho^n_{i_2} (T^{-1(2)}_{(1)})  ...
\rho^n_{i_l} (T^{-1(2)}_{(l-1)}) t^{i_1} t^{i_2} ...t^{i_l}.\nonumber 
\eea 
{From} the definition of $T$ and the relations
$$
(\Delta \otimes \mbox{id})\R=\R_{13}\R_{23}, \qquad
(\mbox{id}\otimes\Delta  )\R=\R_{13}\R_{12}
$$
it follows 
that $T^{-1(1)}\otimes T^{-1(2)}_{(1)}\otimes  ...
\otimes T^{-1(2)}_{(l-1)}$ is a product of $2(l-1)$ 
$\R^{-1}_{mn}$, with suitable $m,n=1,2,...,l$. 
A glance at the explicit form \cite{FadResTak89}  
of the Yang-Baxter matrix $R$ shows that 
$R^{-1}{}^{nn}_{hk}:=\rho^n_h(\R^{-1(1)})\rho^n_k(\R^{-1(2)})
=q^{-1}\delta^n_h\delta^n_k$. 
It follows that  
$$ 
\rho^n_{i_1} (T^{-1(1)}) \rho^n_{i_2} (T^{-1(2)}_{(1)})  ...
\rho^n_{i_l} (T^{-1(2)}_{(l-1)})= q^{-2(l-1)}\delta^n_{i_1}  ...
\delta^n_{i_l}, 
$$ 
which together with the preceding relation gives the recursive relation 
$w_l= q^{3-2l-N}w_{l-1}$; we solve the latter starting 
from $w_1= q^{1-N}$ [see (\ref{wonx})] and we find (\ref{vvalue})$_1$, 
and consequently also (\ref{vvalue})$_2$. 
\ep

\begin{lemma} An element ${\cal O}\in{\cal H}$ is identically zero  
iff for any $f\in \b{R}_q^N$  
\be  
{\cal O}f|=0.  
\ee  
\end{lemma}
{\bf Proof:}  
Let $\{X^{\pi}\}_{\pi\in \Pi}$ be the basis of $\b{R}_q^N$  
dual to the one $\{{\cal D}_{\pi}\}_{\pi\in \Pi}$ of (\ref{Odeco})  
w.r.t. the pairing (\ref{pair}).  
{}From the hypothesis we obtain  
$$  
{\cal O}^{\nu}={\cal O}X^{\nu}|=0\quad \forall \nu\in \Pi\qquad  
\Rightarrow \qquad {\cal O}=\sum\limits_{\nu\in \Pi}{\cal
O}^{\nu}D_{\nu}=0.\qquad\Box   
$$

In order to prove the theorem  we need some more useful relations. 
Let us introduce the short-hand notations 
$$ 
\mu:=1\!+\!q^{2-N},\qquad \qquad\bar\mu:=1\!+\!q^{N-2},\qquad \qquad 
l_z:= \frac{z^l-1}{z-1}, 
$$ 
($l_z$ is called ``$z$-number'' because 
$l_z\stackrel{z\to 1}{\longrightarrow}l$). Moreover, we 
introduce $z$-derivatives (with $z=q,q^{-1}$)
$$
D_zf(r)|:=\frac{f(zr)-f(r)}{(z-1)r}\qquad
\Rightarrow \qquad D_qf(q^{-1}r)|=q^{-1}D_{q^{-1}} f(r)|.
$$
Then, setting henceforth for brevity $\Box:=\partial\cdot\partial$, 
$\hat\Box:=\hat\partial\cdot\hat\partial$,
\bea 
&&\partial^i r^2=\mu\, x^i+q^2r^2 \partial^i,\qquad \qquad  
\partial^i r=\frac{\mu}{1\!+\!q}\frac{x^i}r+q r \partial^i,
\label{dir^2rel}\\
&&\hat\partial^i r^2=\bar\mu\, x^i+q^{-2} r^2
\hat\partial^i,\qquad \qquad  
 \hat\partial^i
r=\frac{\bar\mu}{1\!+\!q^{-1}}\frac{x^i}r+q^{-1} r  
 \hat\partial^i,\nn 
&&\Box x^i=\mu\,\partial^i+q^2 x^i\Box, \qquad \qquad 
\hat\Box x^i=\bar\mu\,\hat\partial^i+q^{-2} x^i\hat\Box, \label{Boxxrel}\\ 
&&\Box\, r^2=\mu^2 \left(q^N\Lambda^{-2}-1\right)(q^2-1)^{-1}
+q^2 r^2\Box, \label{boxr^2rel}\\ 
&&\partial^i f(r)=\frac{\mu}{1\!+\!q}\frac{x^i}r\, D_q f(r)+ 
f(qr)\partial_i, \label{dr^2rel}\\ 
&&\hat\partial^i f(r)=\frac{\bar\mu}{1\!+\!q^{-1}}\frac{x^i}r\,  
D_{q^{-1}} f(r)+ 
f(q^{-1}r)\hat\partial_i.           \label{bardr^2rel} 
\eea 
 
Let 
$$ 
X_l^{i_1...i_l}\!:=r^lS_l^{i_1...i_l} 
={\cal P}^{s,l}{}^{i_1...i_l}_{j_1...j_l} 
x^{j_1}...x^{j_l},\qquad\quad  
T_{l-1}^{i_0i_1...i_l}\! := g^{i_0j_1} 
{\cal P}^{s,l}{}^{i_1...i_l}_{j_1...j_l} 
x^{j_2}...x^{j_l} 
$$ 
[compare with (\ref{defS})]. Clearly  
$r^{1-l}T_{l-1}^{i_0i_1...i_l}\in V_{l-1}$. 
The projector ${\cal P}_{s,l}$ is 
uniquely characterized by the following property \cite{Fio93} 
\be 
{\cal P}_{s,l}{\cal P}_{\pi,m(m\!+\!1)}= 
{\cal P}_{\pi,m(m\!+\!1)}{\cal P}_{s,l}=\delta_{s\pi}{\cal P}_{s,l}, 
\qquad\qquad {\cal P}_{s,l}{}^2={\cal P}_{s,l},
                                               \label{Plproj} 
\ee 
where $\pi=a,s,t$, $m=1,...,l\!-\!1$ and by ${\cal P}_{\pi,m(m\!+\!1)}$  
we have denoted the matrix 
acting as ${\cal P}_{\pi}$ on the $m$-th, $(m\!+\!1)$-th indices and 
as the identity on the remaining ones. Using (\ref{dxrel}) and  
(\ref{Boxxrel}) 
this implies, for $m=1,2,...,l$ 
\bea 
&&{\cal P}^{s,l}{}^{i_1...i_l}_{j_1...j_l}[ 
\partial^{j_m}x^{j_{m\!+\!1}}...x^{j_l}  
-x^{j_m}...x^{j_{l-1}}\partial^{j_l}]=0= 
{\cal P}^{s,l}{}^{i_1...i_l}_{j_1...j_l} 
\partial^{j_m}x^{j_{m\!+\!1}}...x^{j_l}  |\nn  
&&\Box {\cal P}^{s,l}{}^{i_1...i_l}_{j_1...j_l} 
x^{j_m}...x^{j_l} |=0   \nonumber  
\eea 
Using (\ref{dxrel}) (as well as its analog for the $\hat\partial_i$), 
(\ref{gRrel}), (\ref{projectorR}) it follows 
\bea 
&&\partial^{i_0}X_l^{i_1...i_l}| = l_{q^2}T_{l-1}^{i_0i_1...i_l} 
\qquad \qquad\qquad  
x^i\partial_i\,X_l^{i_1...i_l}|=l_{q^2}X_l^{i_1...i_l} \qquad\quad
\label{aa1}\\
&& \hat\partial^{i_0}X_l^{i_1...i_l}| = l_{q^{-2}} 
T_{l-1}^{i_0i_1...i_l} 
\qquad \qquad \qquad\Box X_l^{i_1...i_l}|=0\label{aa2}\\ 
&& x^{i_0}X_l^{i_1...i_l}= 
X_{l+1}^{i_0i_1...i_l}+\frac{r^2 l_{q^2}}{\mu(l-1+N/2)_{q^2}}
T_{l-1}^{i_0i_1...i_l}.   \qquad                  \label{aa3}
\eea 
To prove (\ref{aa3}) note that the decomposition (\ref{ldeco}) 
of the lhs gives (suppressing indices) $x\, X_l=Y_{l+1}+r^2Y_{l-1}$, 
with $Y_j$ combinations of the $X_j$'s. $Y_{l-1}$ can be determined 
applying the Laplacian 
 to both sides and recalling (\ref{aa1}), (\ref{aa2})$_2$: 
\bea 
0&=&\Box(x^{i_0}X_l^{i_1...i_l}-r^2Y_{l-1}^{i_0i_1...i_l})|\nn 
&=&\mu\partial^{i_0}X_l^{i_1...i_l}|- 
\mu^2\left[\left(N/2\right)_{q^2}+ q^N x^i\partial_i\right] 
Y_{l-1}^{i_0i_1...i_l}|\nn 
&=&\mu l_{q^2}T_{l-1}^{i_0i_1...i_l}- 
\mu^2\left[\left(N/2\right)_{q^2}+ q^N(l-1)_{q^2}\right] 
Y_{l-1}^{i_0i_1...i_l}\nn 
&=&\mu \left[l_{q^2}T_{l-1}^{i_0i_1...i_l}- 
\mu(l-1+N/2)_{q^2} 
Y_{l-1}^{i_0i_1...i_l}\right]. \nonumber 
\eea 
Now from (\ref{Plproj}) it follows 
${\cal P}_{s,l\!+\!1}Y_{l-1}\propto{\cal P}_{s,l\!+\!1}T_{l-1}=0$, whence 
$$ 
{\cal P}_{s,l\!+\!1}Y_{l+1}={\cal P}_{s,l\!+\!1}\,x X_l= 
{\cal P}_{s,l\!+\!1}X_{l+1}=X_{l+1}, 
$$ 
and we find that indeed $Y_{l+1}=X_{l+1}$.

\subsubsection*{Proof of Theorem \ref{starsimi}}

Relation (\ref{starsimixi}) is an immediate consequence
of (\ref{starxid})$_1$,  (\ref{barredxi}), (\ref{varphi}). 
The second equality in (\ref{starsimid}) is immediate. As for the first,
\bea
&&\hat\partial^{i_0}\, f(r) X_l^{i_1...i_l}| 
\stackrel{(\ref{bardr^2rel})}{=}
\frac{q\bar\mu}{1\!+\!q}
(D_{q^{-1}}f|)\frac{x^{i_0}}rX_l^{i_1...i_l}
+f(q^{-1}r)\hat\partial^{i_0}\, X_l^{i_1...i_l}| \nn
&& \stackrel{(\ref{aa2}-\ref{aa3})}{=} 
\frac{q^ND_{q^{-1}}f|}{q(1\!+\!q)r}\left[\mu
X_{l+1}^{i_0i_1...i_l}\!+\!\frac{r^2 l_{q^2}}{(l\!-\!1\!+\!\frac N2)_{q^2}}
T_{l-1}^{i_0i_1...i_l}\right]
\!+\!f(q^{-1}r)\,l_{q^{-2}}T_{l-1}^{i_0i_1...i_l}\nonumber
\eea
on one hand, and on the other
\bea
&&v'{}^{-1}\!\partial^{i_0}v'\Lambda\, f(r) X_l^{i_1...i_l}| 
\stackrel{(\ref{Lambdaprop}),(\ref{vvalue})}{=}
v'{}^{-1}\!\partial^{i_0}f(q^{-1}r)X_l^{i_1...i_l}|q^{-(l+N)l/2}\nn
&&\stackrel{(\ref{dr^2rel})}{=}v'{}^{-1}\!\left[q^{-1}\frac{\mu}{1\!+\!q}
(D_{q^{-1}}f|)\frac{x^{i_0}}r+
f(r)\partial^{i_0}\right]X_l^{i_1...i_l}|q^{-(l+N)l/2}\nn
&&\stackrel{(\ref{aa1}-\ref{aa3})}{=}v'{}^{-1}\!\left\{\frac{D_{q^{-1}}f|}{q(1\!+\!q)r}
\left[\mu X_{l+1}^{i_0i_1...i_l}+
\frac{r^2 l_{q^2}}{(l\!-\!1\!+\!\frac N2)_{q^2}}
T_{l-1}^{i_0i_1...i_l}\right]+l_{q^2}
f\,T_{l-1}^{i_0i_1...i_l}\right\}q^{-(l+N)l/2}\nn
&&\stackrel{(\ref{vvalue})}{=}\left\{ \frac{D_{q^{-1}}f|}{q(1\!+\!q)r}\left[\mu 
q^{(l+N-1)(l+1)/2}X_{l+1}^{i_0i_1...i_l}+
\frac{r^2 l_{q^2}\,q^{(l+N-3)(l-1)/2}}{(l\!-\!1\!+\!\frac N2)_{q^2}}
T_{l-1}^{i_0i_1...i_l}\right]\right.\nn
&&+l_{q^2}q^{(l+N-3)(l-1)/2}
f\,T_{l-1}^{i_0i_1...i_l}\Big\}q^{-(l+N)l/2}\nn
&&= \frac{D_{q^{-1}}f|}{(1\!+\!q)r}\left[\mu
q^{(N\!-\!3)/2}X_{l+1}^{i_0i_1...i_l}+
\frac{r^2 l_{q^2}\,q^{-2l-(N-1)/2}}{(l\!-\!1\!+\!\frac N2)_{q^2}}
T_{l-1}^{i_0i_1...i_l}\right]+l_{q^2}q^{(3-N)/2-2l}
f\, T_{l-1}^{i_0i_1...i_l},\nonumber
\eea
whence
\bea
&&[\hat\partial^{i_0}-q^{\frac {N+1}2}v'{}^{-1}\!\partial^{i_0}
v'\Lambda]  f(r) X_l^{i_1...i_l}| =\nn
&&\frac{D_{q^{-1}}f|}{(1\!+\!q)r}\left[\frac{r^2 l_{q^2}\,(q^{N-1}-q^{1-2l})}
{(l\!-\!1\!+\!\frac N2)_{q^2}}T_{l-1}^{i_0i_1...i_l}\right]
+[f(q^{-1}r)-f(r)]\,l_{q^2} q^{2-2l}T_{l-1}^{i_0i_1...i_l}\nn
&&\frac{D_{q^{-1}}f|}{(1\!+\!q)r}\left[r^2 l_{q^2}\,q^{1-2l}(q^2-1)
T_{l-1}^{i_0i_1...i_l}\right]
+(D_{q^{-1}}f|)r(q^{-1}-1)\,l_{q^2} q^{2-2l}T_{l-1}^{i_0i_1...i_l}=0
\nonumber
\eea
leading to $\hat\partial^i=q^{(N\!+\!1)/2} v'{}^{-1}\!\partial^iv'\Lambda$,
equivalent to the claim (\ref{starsimid}). To prove (\ref{starsim}) now 
we just have to proceed as follows. By (\ref{dxirel}) 
$d:=\xi^i\partial^jg_{ij}=q^N\partial^i\xi^jg_{ij}$, whence
\bea
d^{\star} & = & q^N\xi^j{}^{\star}\partial^i{}^{\star}g_{ij}
\stackrel{(\ref{starsimid}),(\ref{starsimixi})}{=}
q^{2N}\xi^kg_{kl}Z'{}^l_j\Lambda^{-2}\Big(-q^{(1-N)/2}v'{}^{-1}\
\partial^jv'\Lambda\Big)\nn
&\stackrel{(\ref{wsigmarel})}{=}& -q^{\frac 32N-\frac
12}\xi^kg_{kl}\Lambda^{-1} q^{1-N}v'{}^{-1} w'\partial^l
w'{}^{-1}v'\nn
&\stackrel{(\ref{eta'drel})}{=}& -q^{\frac 12N+
\frac 12}\xi^kg_{kl} v' q^{-\frac N2-\frac 12}q^{-\eta'{}^2} 
\partial^l q^{\eta'{}^2}   v'{}^{-1} =- \tilde v' \xi^kg_{kl}  
\partial^l\tilde v'{}^{-1} =- \tilde v' d\tilde v'{}^{-1}.\nonumber
\eea

The proof of (\ref{starsimith}) is completely analogous.
\subsection{Proof of formulae (\ref{scalprodp}),(\ref{bla2}) }

\bea
\int_q \alpha_p^{\star}\:{}^*\beta_p | &=&c_k\int_q 
(\theta^{a_1}...\theta^{a_p}\alpha^{\theta}_{a_p...a_1})^{\star}
\theta^{b_{p\!+\!1}}...\theta^{b_N}
\varepsilon_{b_N...b_{p\!+\!1}}{}^{b_1...b_p}\beta^{\theta} _{b_p...b_1}\nn
&=&c_p\int_q \alpha_{a_p...a_1}^{\theta\,\star}\theta^{b_p}...\theta^{b_1}
g_{b_pa_p}...g_{b_1a_1}
\theta^{b_{p\!+\!1}}...\theta^{b_N}
\varepsilon_{b_N...b_{p\!+\!1}}{}^{b_1...b_p}\beta^{\theta} _{b_p...b_1}\nn
&=&c_p\int_q \alpha_{a_p...a_1}^{\theta\,\star}
\varepsilon^{b_p...b_1b_{p\!+\!1}...b_N}
g_{b_pa_p}...g_{b_1a_1}
\varepsilon_{b_N...b_{p\!+\!1}}{}^{b_1...b_p}\beta^{\theta} _{b_p...b_1}\, dV\nn
&=&c_p\int_q \alpha_{a_p...a_1}^{\theta\,\star}
U^{-1}{}^{c_p}_{a_p}...U^{-1}{}^{c_1}_{a_1}
\varepsilon_{c_p...c_1}{}^{b_{p\!+\!1}...b_N}
\varepsilon_{b_N...b_{p\!+\!1}}{}^{b_1...b_p}\beta^{\theta} _{b_p...b_1}\, dV\nn
&=&\frac 1{c_{N\!-\!p}}\int_q \alpha_{a_p...a_1}^{\theta\,\star}
U^{-1}{}^{c_p}_{a_p}...U^{-1}{}^{c_1}_{a_1}
{\cal P}_a{}^{b_1...b_p}_{c_1...c_p}\beta^{\theta} _{b_p...b_1}\, dV\nn
&=&\frac 1{c_{N\!-\!p}}\int_q \alpha_{a_p...a_1}^{\theta\,\star}
U^{-1}{}^{c_p}_{a_p}...U^{-1}{}^{c_1}_{a_1}\beta^{\theta}_{c_p...c_1}\, d^N\!x\nn
&=& \frac{1}{c_{N\!-\!p}}\int_q \alpha^{\theta\,a_p...a_1}{}^{\star} 
\beta^{\theta\,a_p...a_1}{}\,d^N\!x  \qquad\qquad\Box\nonumber
\eea
Here $U$ is the (diagonal, positive-definite) matrix defined
in (\ref{casimirs}).
The second equality is based on the relation \cite{Fio94}
\be
g_{i_1j_1}g_{i_2j_2}...g_{i_Nj_N}\varepsilon^{j_N...j_2j_1}=:
\varepsilon_{i_1i_2...i_N}=\varepsilon^{i_N...i_2i_1}.
\ee

\bea
&&\la {}^*\alpha_p,{}^*\beta_p \ra =
\int_q ({}^*\alpha)^{\star}\:{}^*{}^*\beta_p\nn
 &&=c_p^{\star}\int_q 
(\alpha^{\theta}_{a_p...a_1})^{\star}
\varepsilon_{a_N...a_{p\!+\!1}}{}^{a_1...a_p}
\theta^{b_N}...\theta^{b_{p\!+\!1}}g_{b_{p\!+\!1}a_{p\!+\!1}}
...g_{b_Na_N}\theta^{b_1}...\theta^{b_p}\beta^{\theta}_{b_p...b_1}\nn
 &&=c_p\,
\varepsilon_{a_N...a_{p\!+\!1}}{}^{a_1...a_p}
\varepsilon^{b_N...b_{p\!+\!1}b_1...b_p}
\int_q (\alpha^{\theta}_{a_p...a_1})^{\star}g_{b_{p\!+\!1}a_{p\!+\!1}}
...g_{b_Na_N}\beta^{\theta}_{b_p...b_1}dV\nn
 &&=c_p\,
\varepsilon_{a_N...a_{p\!+\!1}}{}^{a_1...a_p}
\varepsilon^{b_N...b_{p\!+\!1}b_1...b_p}
\int_q (\alpha^{\theta}_{a_p...a_1})^{\star}g_{b_{p\!+\!1}a_{p\!+\!1}}
...g_{b_Na_N}\beta^{\theta}_{b_p...b_1}dV\nn
&& \qquad \qquad g_{c_1d_1}...g_{c_pd_p}
g^{d_1b_1}...g^{d_pb_p} \nn
 &&=c_p\,
\varepsilon_{a_N...a_{p\!+\!1}}{}^{a_1...a_p}
\varepsilon^{d_p...d_1a_{p\!+\!1}...a_N}
\int_q (\alpha^{\theta}_{a_p...a_1})^{\star}g^{d_1b_1}...g^{d_pb_p} 
\beta^{\theta}_{b_p...b_1}dV\nn
 &&=c_p\,
\varepsilon_{d_p...d_1}{}^{a_{p\!+\!1}...a_N}
\varepsilon_{a_N...a_{p\!+\!1}}{}^{a_1...a_p}
U^{-1}{}^{d_p}_{b_p}...U^{-1}{}^{d_1}_{b_1}
\int_q (\alpha^{\theta}_{a_p...a_1})^{\star}\beta^{\theta}_{b_p...b_1} dV\nn
 &&=\frac 1{c_{N\!-\!p}}
{\cal P}_a{}^{a_1...a_p}_{d_1...d_p}
U^{-1}{}^{d_p}_{b_p}...U^{-1}{}^{d_1}_{b_1}
\int_q (\alpha^{\theta}_{a_p...a_1})^{\star}\beta^{\theta}_{b_p...b_1} dV
\nn &&=\frac 1{c_{N\!-\!p}} 
U^{-1}{}^{d_p}_{b_p}...U^{-1}{}^{d_1}_{b_1} 
\int_q (\alpha^{\theta}_{d_p...d_1})^{\star}\beta^{\theta}_{b_p...b_1} dV 
=rhs(\ref{bla2})
\qquad\quad\Box\nonumber 
\eea

\subsection{Studying the positivity relation (\ref{positi})}

According to odd or even $p=N/{\gamma}$ the matrix elements of $M(a)$ will take the two
different forms
\bea
&&M^j_{j'} \!=\!\frac {\pi^2}{8h}\!\!\int\limits^{\infty}_{-\infty}\!\! 
d\omega'{}  \int\limits^{\infty}_{-\infty}\!\! 
d\omega \int\limits^{\infty}_{-\infty} \!\!  dy \, \check{m}(y)
\frac{e^{i\frac{\pi}2[(\omega'\!-\!\omega)y+\omega j\!-\!\omega'j']
-\frac {a}h\omega'{}^2\!+\!\frac
N2h(j\!+\!j'\!+\!2\beta)}}{C(\frac{\pi^2}{\gamma
h}\omega)C(\frac{\pi^2}{h\gamma}\omega')} \nn
 && C(\omega):=\cases{\cosh(\omega)
 \quad\mbox{ if }p:=N/{\gamma}\mbox{ isodd}\cr
 \sinh(\omega) \quad\mbox{ if }p:=N/{\gamma}\mbox{ is even.}}
\label{I4} 
 \eea
To obtain  the previous formula from (\ref{I3})
we have also performed the change of integration variables
$y\to h(y/2\!+\!\beta)$, $\omega \to\pi\omega/h$, $\omega' \to\pi\omega'/h$
and set $\check{m}(y):=\tilde m(hy/2\!+\!h\beta)$, whence it follows
for any $k\in\b{Z}$  
$$
\check m(y+2k)=\check m(y),\qquad\qquad\check m(-y)=\check m(y),
$$
so that
$$
\check m(y)=\sum_{k=-\infty}^{\infty}m_ke^{ik\pi y},
\qquad\qquad\mbox{with }m_{-k}=m_k=m_k^{\star}. 
$$
We also  define
$$
\check{\phi}(\omega):=\sum\limits_{j\in\b{Z}}e^{-i\pi\omega j\!+\!\frac
N2h(j\!+\!\beta)} R^j_{\phi} \qquad\Rightarrow \qquad
\check{\phi}(\omega+2k)=\check{\phi}(\omega), \quad \forall k\in\b{Z}.
$$
Replacing in (\ref{phipsi'}) (with $\psi=\phi$) we find
\bea
\la \phi, q^{a\eta'{}^2}\phi\ra'&= &\frac {\pi^2}{8h}
\!\!\int\limits^{\infty}_{-\infty}\!\! 
d\omega'{}  \int\limits^{\infty}_{-\infty}\!\! 
d\omega \int\limits^{\infty}_{-\infty} \!\!  dy \, \check{m}(y)
\frac{e^{i\frac{\pi}2[(\omega'\!-\!\omega)y]
-\frac {a\pi^2}h\omega'{}^2} [\check{\phi}(\omega)]^{\star}
\check{\phi}(\omega')}{C(\frac{\pi^2}{h\gamma}\omega)C(\frac{\pi^2}{h\gamma}\omega')} \nn
&= &\frac {\pi^2}{2h}
\!\!\int\limits^{\infty}_{-\infty}\!\! 
d\omega'{}  \int\limits^{\infty}_{-\infty}\!\! 
d\omega 
\sum_{k=-\infty}^{\infty}m_k\,
\delta(\omega'\!-\!\omega\!+\!2k)
\frac{e^{-\frac {a\pi^2}h\omega'{}^2} [\check{\phi}(\omega)]^{\star}
\check{\phi}(\omega')}{C(\frac{\pi^2}{h\gamma}\omega)C(\frac{\pi^2}{h\gamma}\omega')} \nn
&= &
\frac {\pi^2}{2h}\sum_{k=-\infty}^{\infty}m_k
\int\limits^{\infty}_{-\infty}\!\! 
d\omega'{} \frac{e^{-\frac {a\pi^2}h\omega'{}^2}
|\check{\phi}(\omega')|^2}{C[\frac{\pi^2}{h\gamma}(\omega'\!+\!2k)]C(\frac{\pi^2}{h\gamma}\omega')} 
\nn &= &
\frac {\pi^2}{4h}\int\limits^{1}_{-1} \!\!  dy \,
\check{m}(y) \int\limits^{\infty}_{-\infty}\!\! 
d\omega{} \frac{|\check{\phi}(\omega)|^2e^{-\frac {a\pi^2}h\omega{}^2}}{C(\frac{\pi^2}{h\gamma}\omega)} 
\sum_{k=-\infty}^{\infty}\frac{\cos(k\pi
y)}{C[\frac{\pi^2}{h\gamma}(\omega\!+\!2k)]} \nn &= &
\frac {\pi^2}{4h}\int\limits^{1}_{-1}\!\! 
d\omega |\check{\phi}(\omega)|^2\int\limits^{1}_{-1} \!\!  dy \,
\check{m}(y) \, K\left(\omega,y,\frac{\pi^2}{h\gamma},\frac
{a\pi^2}h \right).
\eea
Thus $\la \phi, q^{a\eta'{}^2}\phi\ra'$ will be positive for any $\phi$
if
\be
\int\limits^{1}_{-1} \!\!  dy \,
\check{m}(y) \, K\left(\omega,y,\frac{\pi^2}{h\gamma},\frac
{a\pi^2}h \right)>0 \qquad\qquad \forall \omega\in]-1,1], \label{posicon}
\ee
where
$$
K(\omega,y, \delta,t):= 
\sum_{l=-\infty}^{\infty}
\frac{e^{-t(\omega\!+\!2 l)^2}}{C[\delta(\omega\!+\!2l)]}  
 \sum_{k=-\infty}^{\infty}
\frac{\cos(k\pi y)}{C[\delta(\omega\!+\!2(k\!+\!l)]}  .
$$

The weight characterizing the Jackson integral, 
$\check{m}(y)\sim\sum_{l=-\infty}^{\infty}\delta(y-2l)$ certainly fulfills
(\ref{posicon}) for any choice of $a,h,\gamma$ because the integral appearing
there  reduces to $K(y=0)$ which is
manifestly
 positive. In fact, by continuity, $K$ will remainpositive
 at least in a neighbourhood of $y=0$, so that (\ref{posicon}) willbe
 fulfilled also by weights $\check{m}(y)$ nonvanishing on some 
suitable interval including $y=0$. A more detailed
characterization of weights  $\check{m}(y)$ and parameters $a,h,\gamma$ 
such that (\ref{posicon}) is fulfilled is left as a
possible subject for future work. If $K$ were strictly positive for all $y$
all weights  $\check{m}(y)$ would do the job. 
Note also that
for $h\to 0$ one finds
$\sqrt{a/h\pi} e^{-\frac {a\pi^2}h\omega{}^2}\sim \delta(\omega)$
and also
$1/C(k\pi^2/{\gamma}h)\to \delta_k^0$  and therefore
$$
\la \phi, q^{a\eta'{}^2}\phi\ra'\sim m_0 \frac{\pi^{\frac 32}}{2\sqrt{ah}}
\left|\check{\phi}(0)\right|^2\ge 0
$$
which is nonnegative for any $\phi$ and any choice of $\check{m}(y)$.

\end{document}